%% file: main.tex
\documentclass{article}

\usepackage[english]{babel}
\usepackage{parskip}
\usepackage[a4paper, total={6in, 8in}]{geometry}

\usepackage{amsmath}
\usepackage{amsfonts}
\usepackage{amssymb}
\usepackage{xcolor}
\usepackage{amsthm}
\usepackage{graphicx}
\usepackage{bm}
\usepackage{newclude}
\usepackage[colorlinks=true, allcolors=blue]{hyperref}
\usepackage[toc,page]{appendix}
\usepackage{mathtools}
\usepackage{comment}
\usepackage{authblk}
\usepackage{enumerate}

\usepackage{pgfplots}
\usepackage{pgfplotstable}
\pgfplotsset{compat=1.7}
\usepackage{xcolor}
\usepackage{wasysym}
\usepackage{customcolors}

\usepackage{algorithm}
\usepackage{algpseudocode}
\usepackage{caption}
\usepackage{subcaption}
\usepackage{comment}
\usepackage{afterpage}

\newcommand{\anorm}[1]{\|#1\|_{{V}}}
\newcommand{\amin}{a_{min}}
\newcommand{\ch}[2]{\binom{#1}{#2}}
\newcommand{\Avar}[1]{\partial^{\mu-\nu}{A}(\bm{#1})}

\newcommand{\A}{\Avar{0}}

\newcommand{\abs}[1]{\left|#1\right|}
\newcommand{\x}{{{\bm{x}}}}
\newcommand{\xx}{{x}_1}
\newcommand{\xy}{{x}_2}
\renewcommand{\r}{|\x|}
\renewcommand{\a}{\frac{r_0^-}{4}}
\newcommand{\nui}{{\nu^{(1)}}}
\newcommand{\nuii}{{\nu^{(2)}}}
\newcommand{\nuiii}{{\nu^{(3)}}}
\newcommand{\nuiv}{{\nu^{(4)}}}

\newcommand{\Ca}{{{C}_a}}
\newcommand{\Cf}{{{C}_f}}
\newcommand{\Cat}{{\tilde{C}_a}}
\newcommand{\Cft}{{\tilde{C}_f}}

\newcommand{\Cpc}{C_p}
\newcommand{\Csigma}{C_{\sigma}}

\newcommand{\Kb}{\Kt}

\newcommand{\Kt}{K_T}
\newcommand{\uparamdomain}{\mathrm{u}}

\newcommand{\chibar}{\bar{\chi}}
\newcommand{\chibarone}{\chibar_1}
\newcommand{\chibartwo}{\chibar_2}

\newtheoremstyle{break}
  {.5em}
  {}
  {\itshape}
  {}
  {\bfseries}
  {.}
  {.5em}
  {}%

\theoremstyle{break}
\newtheorem{theorem}{Theorem}[section]
\newtheorem{lemma}[theorem]{Lemma}
\newtheorem{corollary}[theorem]{Corollary}
\newtheorem{definition}[theorem]{Definition}
\newtheorem{assumpttion}[theorem]{Assumption}
\newtheorem{example}[theorem]{Example}
\newtheorem{remark}[theorem]{Remark}

\numberwithin{equation}{section}
\numberwithin{figure}{section}

\DeclareMathOperator{\tr}{tr}
\DeclareMathOperator{\sign}{sign}
\renewcommand{\d}{\mathrm{d}}
\DeclareMathOperator{\supp}{supp}
\DeclareMathOperator*{\argmax}{arg\,max}

\newcommand{\dx}{\mathop{}\!\mathrm{d} {\bm{x}}}
\renewcommand{\d}{\mathop{}\!\mathrm{d}}
\newcommand{\D}{\mathop{}\!\mathrm{D}_x}

\newcommand{\Dy}{\mathop{}\mathcal{D}(\y)}
\newcommand{\y}{\mathop{}\!\bm{y}}

\newcommand{\deltatilde}{{\tilde{\delta}}}

\allowdisplaybreaks

\usepackage{mathtools}

\mathtoolsset{showonlyrefs}

\newif\ifgraphs
\graphsfalse

\title{Exploiting locality in sparse polynomial approximation of parametric elliptic PDEs and application to parameterized domains}

\author{Wouter van Harten, Laura Scarabosio}

\affil{\small Institute for Mathematics, Astrophysics and Particle Physics, Radboud University, Nijmegen, The Netherlands}

\providecommand{\keywords}[1]{\textbf{Keywords}: {#1}}
\providecommand{\ams}[1]{\textbf{AMS subject classification}: {#1}}
\begin{document}

\maketitle

\begin{abstract}
This work studies how the choice of the representation for parametric, spatially distributed inputs to elliptic partial differential equations (PDEs) affects the efficiency of a polynomial surrogate, based on Taylor expansion, for the parameter-to-solution map. In particular, we show potential advantages of representations using functions with localized supports. As model problem, we consider the steady-state diffusion equation, where the diffusion coefficient and right-hand side depend smoothly but potentially in a \textsl{highly nonlinear} way on a parameter $\y\in [-1,1]^{\mathbb{N}}$. Following previous work for affine parameter dependence and for the lognormal case, we use pointwise instead of norm-wise bounds to prove $\ell^p$-summability of the Taylor coefficients of the solution. 
As application, we consider surrogates for solutions to elliptic PDEs on parametric domains. Using a mapping to a nominal configuration, this case fits in the general framework, and higher convergence rates can be attained when modeling the parametric boundary via spatially localized functions. The theoretical results are supported by numerical experiments for the parametric domain problem, illustrating the efficiency of the proposed approach and providing further insight on numerical aspects. Although the methods and ideas are carried out for the steady-state diffusion equation, they extend easily to other elliptic and parabolic PDEs.

\end{abstract}


\ams{35B30, 41A10, 41A58, 41A63, 65N15, 65T60.}
\include*{1introduction/1introduction}

\include*{2summability/2summability}

\include*{4model_problem/4model_problem}

\include*{4application/4application}

\include*{5numerical_results/5numerical_results}

\include*{6conclusion/6conclusion}

\bibliographystyle{plain} 
\bibliography{bibtex/library} 

\newpage
\include*{7appendices/7appendices}

\end{document}

%% file: 1introduction/1introduction.tex
\section{Introduction}
\label{sec:introduction}
Our growing demand for more efficient and accurate devices induces a need for computational prototyping in multi-query \cite{Rozza2008} or outer-loop \cite{Peherstorfer2018} scenarios, such as optimization, inference, and uncertainty quantification. These tasks often involve the evaluation of a mathematical model for many values of an input parameter, for instance, the quantity to be optimized or an uncertain quantity. The model evaluation is often the computational bottleneck in these frameworks, and \textsl{surrogate models} are a key tool to reduce the computational cost dramatically.

In this paper, we focus on scenarios where the system's behavior can be modeled by a \textsl{partial differential equation} (PDE) of \textsl{elliptic} type with spatially distributed, parametric inputs. More precisely, we consider the parametric model problem
\begin{equation}
	\begin{cases}
		-\nabla\cdot\left(A(\y)\nabla u(\y) \right) = F(\y) & \text{in } D,\\
		u = 0 & \text{on } \partial D,\\
		\text{for every }\bm{y} \in Y := [-1, 1]^\mathbb{N},
	\end{cases} \label{eq:pde}
\end{equation}
where $\y$ is the \textsl{high-dimensional} parameter and $D\subset \mathbb{R}^d$ is a bounded, connected domain for $d\in\mathbb{N}$ (usually $d=1,2,3$ in applications) with Lipschitz boundary. We require that, for every $\y\in Y$, $A(\y)\in\left(L^{\infty}(D)\right)^{d\times d}$ and $F(\y)\in L^2(D)$, with a parameter-independent bound on their norms. Furthermore, $A(\y,\cdot)$ is positive definite, uniformly with respect to the spatial coordinate and the parameter. 

The parameter dependence of the input, entering \eqref{eq:pde} via the PDE coefficient and right-hand side, propagates to the solution $u$, with $u(\y)\in H_0^1(D)$ for every $\y\in Y$, and possible output functionals depending on it. We are interested in constructing efficiently a surrogate, that is an approximant, for the parameter-to-solution map $\y \mapsto u(\y)$. We focus on sparse \textsl{polynomial} surrogates, and in Section \ref{sec:conclusion} we highlight potential implications for approximation properties of other surrogates and quadrature. A consistent part of this work addresses the application to elliptic PDEs on \textsl{parameter-dependent domains} (with fixed coefficients and right-hand side), which, as we will see, can be recast in the setting \eqref{eq:pde}. Motivated by the application to parametric domains, among many others, we are particularly interested in the case that $A$ and $F$ in \eqref{eq:pde} depend smoothly but possibly in a highly nonlinear way on the parameter $\y$.

The first results on convergence rates for polynomial surrogates of the parameter-to-solution map considered the case of a finite-dimensional parameter, see for instance \cite{Babuska2007,Babuska2004,Xiu2006a}. Convergence results for a parameter with a countable number of entries, in other words, dimension-independent rates, were first obtained in \cite{Todor2007} and improved in \cite{Cohen2010,Cohen2011}. These works study the model problem \eqref{eq:pde} when the right-hand side is independent of the parameter and the coefficient depends on it in an affine way, that is $A(\y)=a(\y)I_d$, where $I_d$ is the identity matrix and $a(\y)=a_0 + \sum_{j\geq 1}y_j\psi_j$, for given $a_0, \left\{\psi_j\right\}_{j\geq 1}\in L^{\infty}(D)$ such that $a$ is uniformly positive and bounded (extension to matrix-valued coefficients is straightforward). From estimates on Taylor coefficients (with respect to the parameter), the aforementioned works obtain convergence rates on polynomial chaos and collocation-based surrogates. 

A significant step forward was taken in \cite{Chkifa2015}, where convergence rates for Legendre instead of Taylor coefficients were proved to hold for a much larger class of problems, namely \eqref{eq:pde} with smooth but potentially highly nonlinear dependence of $A(\y)$ and $F(\y)$ on the parameter, and a broader class of PDEs. We refer to \cite{Cohen2015} and \cite[Ch. 2-4]{Adcock2022} for a review of the so-far mentioned results. 

All these works are based on the fact that the solution to the PDE, as a function of the parameter, admits a smooth extension on sets of the form $\mathcal{O}=\otimes_{j\geq 1}\mathcal{O}_j$ such that $[-1,1]\subset\mathcal{O}_j$ and $\mathcal{O}_j$ is an open set on the real line or in the complex plane, and the sizes of the sets $\mathcal{O}_j$ are determined by the sensitivity of the parametric input with respect to $y_j$, when sensitivity is measured using Banach space \textsl{norms}. More concretely, if the coefficient and right-hand side in \eqref{eq:pde} depend on a parametric input $\mathfrak{p}(\y) = \mathfrak{p}_0 + \sum_{j\geq 1}y_j\psi_j$ taking values in a Banach space $X$, where $\mathfrak{p}(\y)=A(\y)$ in the affine case but could also be something more involved such as the boundary of a domain \cite{Castrillon-Candas2016,Harbrecht2016,Hiptmair2018}, then, in the works mentioned so far, the size of $\mathcal{O}_j$ depends on the decay of $\lVert\psi_j\rVert_X$ as $j\rightarrow\infty$. In particular, they do not take into account the support of $\psi_j$, $j\geq 1$. Results in the same spirit, although requiring quite different proof techniques, hold in the lognormal case, that is when $A(\y)$ in \eqref{eq:pde} is the exponential of a Gaussian random field, and the parameter $\y$, taking values in the whole $\mathbb{R}^{\mathbb{N}}$ with underlying Gaussian measure, corresponds to the image of the Gaussian random variables in the field expansion, see for instance \cite{Ernst2018,Hoang2014}. 

In the last years, starting from \cite{Bachmayr2017a}, it has become clear that using \textsl{pointwise} instead of norm-wise bounds on the decay of $\psi_j$ as $j\rightarrow\infty$ can lead to a faster decay of the Taylor or Legendre coefficients of the PDE solution when the supports of $\psi_j$, $j\geq 1$, are localized. The main idea is that, in this case, using norms to measure the effect of each parameter coordinate on the variability of the output might be too pessimistic. One should exploit instead that, at a specific location in the domain $D$, only few functions contribute significantly to changes of the input $\mathfrak{p}$ (and so, in some way, to changes in the PDE solution).  The work \cite{Bachmayr2017a} considers the model problem \eqref{eq:pde} and an \textsl{affine} dependence of the diffusion coefficient on the parameter (and fixed right-hand side). Similar results have been shown in the lognormal case \cite{Bachmayr2017,Dung2022}. 

This paper builds upon \cite{Bachmayr2017,Bachmayr2017a}. Our main contributions are: (a) we extend the results in \cite{Bachmayr2017a}, which are for the affine case, to the case that the coefficient and right-hand side in \eqref{eq:pde} depend \textsl{nonlinearly}, but smoothly, on the parameter; (b) we show how our results can be applied to treat elliptic problems on parametric domains, and in particular advantages, in terms of convergence rates for surrogates, of using functions with localized supports to represent the parametric boundary, or more generally domain variations; (c) for the application to parametric domains, we present numerical results that, apart from illustrating the theory, provide further insight on computational approaches to handle the parameter-dependent geometry.

In the application to parametric domains, we adopt the mapping approach developed within the framework of uncertainty quantification by Tartakovsky and Xiu \cite{Tartakovsky2006, Xiu2006}. This involves transforming the PDE on the parameterized domain to a parameterized PDE on a nominal (reference) domain. Such technique is often used in shape optimization \cite{Onyshkevych2021,Haubner2020,Hiptmair2015,Onyshkevych2021a} and provides a framework for which theory and discretization algorithms are well established. 
In particular, when applied to a stationary diffusion equation, the resulting PDE on the nominal configuration fits into the framework \eqref{eq:pde}. The smoothness of the solution on the reference configuration with respect to the parameter has been analyzed in \cite{Castrillon-Candas2016,Harbrecht2016,Multerer2019} for the linear elliptic equation, in \cite{Hiptmair2018} for a low-frequency Helmholtz transmission problem, in \cite{Castrillon-Candas2021} for the linear parabolic equation, in \cite{Cohen2018} for the stationary Navier-Stokes equations, and in \cite{Jerez2017} for Maxwell's equations in frequency domain. 

While all these works prove smoothness in terms of the global properties of the parametric boundary, in the spirit of the general results in \cite{Babuska2007, Chkifa2015}, the present work relies on \textsl{pointwise} properties, in the spirit of \cite{Bachmayr2017,Bachmayr2017a,Dung2022}, showing that parametrizing the boundary, or, more generally, the domain mapping, using functions with localized supports leads to higher convergence rates for surrogates based on truncated Taylor expansions. Our proofs are based on what in \cite{Dung2022} is referred to as ``bootstrap" arguments, that is iterated differentiation, which, for the application to parametric domains, are in the spirit of the proofs in \cite{Harbrecht2016}, although using different estimates. We highlight that, in the case of parametric interface problems, one might also be interested in the solution on the physical, not the nominal, configuration. In that case, the dependence of the solution on the parameter is not smooth \cite{Motamed2013,Scarabosio2017,Scarabosio2022}. 

The literature mentioned so far is for variational formulations on the domain, for analyses with boundary integral equations we refer to \cite{Dolz2023b,Henriquez2021,Pinto2023}. 
Other approaches different from the mapping are possible to handle parameter-dependent geometries. The fictitious domain approach \cite{Canuto2007} and level set methods \cite{Nouy2008,Nouy2007, Osher2001} embed the parameter-dependent domain in a larger, hold-all domain, and level set methods can handle both geometrical and topological changes. However, since these techniques recast a problem on a parameter-dependent domain to a problem with a parametric interface, the issue raised in \cite{Motamed2013,Scarabosio2017,Scarabosio2022} appears. Namely, the solution depends non-smoothly on the parameter in a neighborhood of the boundary, affecting the convergence of algorithms for polynomial-based surrogates. In uncertainty quantification, when one is not interested in the surrogate, but in moments of the quantities of interest only, and when the shape variations are small, it is possible to adopt perturbation methods \cite{Mathematics2017,Harbrecht2013,Harbrecht2008}.

When discussing the smoothness of the solution on the nominal configuration, we analyze theoretically and compare numerically two approaches to construct a domain mapping starting from the boundary parametrization, namely an explicit expression using a mollifier and the harmonic extension \cite{Li2001,Xiu2006}. The latter shows better localization properties in the numerical results, and it is better suited for complex geometries. We analyze it as a prototype for other PDE-based methods, such as those based on elasticity equations, see for instance \cite{Cizmas2008,Dwight2009}. Another possibility, acting on the discretized PDE, is to use isogeometric analysis (IGA), possibly with boundary elements \cite{Dolz2022,Dolz2023}.

This work starts with the introduction of Taylor approximations to the parameter-to-solution map for problem \eqref{eq:pde} in Section \ref{sec:summability}. In the remainder of the latter, we prove  $\ell^p$-summability of the Taylor coefficients via a weighted $\ell^2$-summability argument. In Section \ref{sec:modelproblem}, we introduce the stationary diffusion equation on parameterized domains, and we use the results from the previous section to investigate the summability properties of the Taylor coefficients for the solution on the nominal configuration. There, we also discuss different mapping approaches. In Section \ref{sec:numerical}, we finish the core of our paper with numerical experiments illustrating our theoretical results for the parameterized domain problem. We end with some final remarks and conclusions in Section \ref{sec:conclusion}.

Throughout the paper, we use the following conventions. We denote the spectral norm of a $n$-tensor by $\left|\,\cdot\,\right|_{n,2}$, and, for the special case $n=1$, we use the short notation $\left|\,\cdot\,\right|$. For parameter- and space-dependent quantities, when using the aforementioned norms in proofs, in order to lighten the notation we will omit explicit dependence on the space coordinate. For instance, for the diffusion coefficient in \eqref{eq:pde}, $\left|A(\y)\right|_{2,2}$ will denote the spectral norm of the matrix $A(\y,\x)$, for given $\y\in Y$ and $\x\in D$. For the norm on an infinite-dimensional Banach space $X$, instead, we use $\lVert\,\cdot\,\rVert_X$.

%% file: 2summability/2summability.tex
\section{Summability of Taylor coefficients}
\label{sec:summability}
We are concerned with problem \eqref{eq:pde} when the solution map $\y \mapsto u(\y)$ from $Y$ to $V:=H_0^1(D)$ is analytic, as better specified below. In this case, we can introduce a polynomial expansion mapping $\y$ to $u(\y)$, and in particular we can use the Taylor series
\begin{equation}
	{u}(\y) = \sum_{\mu\in\mathcal{F}} t_\mu \y^\mu, \label{eq:tayloreq}
\end{equation}
where $t_\mu \in V$ are the Taylor coefficients, defined by 
\begin{equation}
	t_\mu = \frac{1}{\mu!}\partial^\mu {u}(\bm{0}),\label{eq:taylordef}
\end{equation}
and $\mathcal{F}$ is the set of all finitely supported multi-indices:
\begin{align*}
	\mathcal{F}=\{\mu=(\mu_1,\mu_2,\mu_3,...), \mu_i \in \mathbb{N}_0, \text{supp}(\mu) < \infty \},
\end{align*} 
with $\mathbb{N}_0:=\mathbb{N}\cup \left\{0\right\}$. The support of a multi-index $\mu$ is defined by $\text{supp}(\mu) := \max\{k \in \mathbb{N}: \mu_k > 0\}$. In \eqref{eq:tayloreq}, $\y^\mu = \prod_{j\geq 1}y_j^{\mu_j}$, and in \eqref{eq:taylordef} $\mu! = \prod_{j \geq 1} \mu_j!$. Other mathematical operations can be extended to multi-indices in a similar manner \cite{Cohen2010}. For example, we have the binomial coefficient $\binom{\mu }{ \nu} = \prod_{j \geq 1}\binom{\mu_j }{ \nu_j}$; ordering $\mu \leq \nu \iff \mu_j \leq  \nu_j $ for all $ j \in \mathbb{N}$; strict ordering $\mu < \nu \iff \mu \leq  \nu$ and $\mu \neq \nu$; and the order $|\mu| = \sum_{j \geq 1} \mu_j$.

A surrogate for the parameter-to-solution map can be obtained by selecting, in \eqref{eq:tayloreq}, a finite-dimensional subset of $\mathcal{F}$ corresponding to the Taylor coefficients with largest $V$-norm. The convergence rate of such surrogate can be derived using Stechkin's lemma, as we will also recall in Corollary \ref{cor:stechkin}, which relies on the $\ell^p$-summability of $(\|t_\nu\|_{V})_{\nu \in \mathcal{F}}$. 

Following \cite{Bachmayr2017a}, to show $\ell^p$-summability of 
\begin{align}
	(\|t_\mu\|_{V})_{\mu \in \mathcal{F}}, \label{eq:taylorseq}
\end{align}
we will first formulate a variational formulation for the Taylor coefficients in Section \ref{sec:taylor}. 
With this formulation at hand, in Section \ref{sec:weightedl2}, we will show the  summability of 
\begin{equation}
	(\rho^{2\mu}\|t_\mu\|_{V}^2)_{\mu\in \mathcal{F}}, \label{eq:weightedtaylorseq}
\end{equation}
 with a positive sequence $\rho=(\rho_j)_{j \geq 1}$ and $\rho^{2\mu}$ defined elementwise, similarly to $\y^\mu$. From the weighted $\ell^2$-summability, we will conclude the $\ell^p$-summability of $(\|t_\mu\|_{V})_{\mu\in \mathcal{F}}$ in Section \ref{sec:ellpsum}.

\subsection{Variational formulation for the Taylor coefficients}
\label{sec:taylor}
Before we can derive a variational formulation for the derivatives of ${u}(\y)$ with respect to $\y$, we need to ask for analyticity and boundedness of the diffusion coefficient and the right-hand side of equation \eqref{eq:pde} in the following assumption:
\begin{assumpttion} \label{ass:afanalytic}
	The functions $A(\y;\bm{x})$ and $F(\y;\bm{x})$ in \eqref{eq:pde} are such that:
\begin{enumerate}[(a)]
		\item the quantities $\|A(\y, \cdot)\|_{L^\infty(D)}$ and $\|F(\y, \cdot)\|_{L^2(D)}$ have $\y$-independent upper bounds, for $\y \in Y$;
		\item they are analytic as maps from $\mathcal{O}:=\bigotimes_{j\geq 1}\mathcal{O}_j$ to $\left(L^\infty(D)\right)^{d\times d}$ and $L^2(D)$, respectively, where each $\mathcal{O}_j\subset\mathbb{R}$ is an open interval containing $[-1,1]$;
		\item there exists a constant $A_{min}>0$ such that, for a.e. $\bm{x} \in D$ and all $\y\in \mathcal{O}$, it is an lower bound of the smallest singular value of $A(\y; \x)$.
\end{enumerate}
\end{assumpttion}
Assumption \ref{ass:afanalytic}(a) guarantees   that $u \in L^\infty(Y, V)$. Now, we derive the variational formulation for the derivatives $\partial^\mu {u}(\y)$. 
First, we obtain an important lemma, giving us a variational formulation for the Taylor coefficients in \eqref{eq:tayloreq}: 
\begin{lemma}
	\label{col:taylorvar}
	Let Assumption \ref{ass:afanalytic} hold. Then, for every $\mu \in \mathcal{F}$, the Taylor coefficient $t_\mu \in V$ defined in \eqref{eq:taylordef} is the unique solution to  
	\begin{align}
		\int_{{D}}  A(\bm{0}) \nabla t_{\mu} \cdot \nabla v \dx = -\sum_{\nu\in S_\mu} \frac{1}{(\mu - \nu)!} \int_{{D}} \partial^{\mu - \nu} {A}(\bm{0}) \nabla t_\nu  \cdot \nabla v \dx + \frac{1}{\mu!} \int_{{D}} &\partial^\mu {F}(\bm{0})  v \dx, \nonumber \\
		&\text{ for all } v \in V, \label{eq:tayorweakform}
	\end{align}
	with $S_\mu = \{\nu \in \mathcal{F} : \nu < \mu\}$. 
\end{lemma}
\begin{proof}
Thanks to Assumption \ref{ass:afanalytic}, for $e_j=(0, \cdots, 0, 1, 0, \cdots )\in\mathcal{F}$ and every $\y \in Y$ the quotient $\frac{u(\y+he_j)-u(\y)}{h}\in V$ is well-defined, for $|h|$ small enough. By argumentations similar to those in \cite{Cohen2010}, passing to the limit $|h|\rightarrow 0$, one can show that the derivative $\partial^{e_j} u(\y) \in V$ is well-defined as the unique solution to the variational formulation 
	\begin{align}
	    &\int_{{D}} {A}(\bm{y}) \nabla \partial^{e_j}{u}(\bm{y}) \cdot\nabla v \dx = L_0(v), \quad \text{ for all } v \in V,\label{eq:singleder}
	\end{align}
	with 
	\begin{align}
	    L_0 (v)= -\int_{{D}} \partial^{e_j}{A}(\bm{y}) \nabla {u}(\bm{y})\cdot\nabla v\dx + \int_{{D}}\partial^{e_j} {F}(\bm{y}) v\dx.
	\end{align}
	By repeated application, we arrive at the variational formulation for the general derivative:
	\begin{align}
    	\int_{{D}}  {A}(\bm{y}) \nabla\partial^{\mu} u(\bm{y})\cdot \nabla v \dx = L(v) ,\nonumber \quad \text{ for all } v \in V, \\
    	L(v) =  -\sum_{\nu\in S_\mu} \ch{\mu}{\nu} \int_{{D}} \partial^{\mu - \nu} {A}(\bm{y}) \nabla \partial^{\nu} u(\bm{y}) \cdot \nabla v \dx + \int_{{D}}  &\partial^\mu {F}(\bm{y})v \dx. \label{eq:gender}
	\end{align}
	Equation \eqref{eq:tayorweakform} is then obtained by evaluating the formulation above at $\y=\bm{0}$ and dividing both sides by $\mu!$. 
\end{proof}
With these results in our toolbox, we are interested in the convergence properties of the Taylor coefficients \eqref{eq:taylorseq}.
The goal of the next subsection is to show that we have summability of the weighted sequence \eqref{eq:weightedtaylorseq}.

\subsection{Weighted summability of Taylor coefficients}
\label{sec:weightedl2}
The weighted $\ell^2$-summability of the Taylor coefficients is given by Theorem \ref{thm:l2summability}, which we will prove in the rest of this section. We note that, differently from Lemma \ref{col:taylorvar}, where we were only concerned with the existence of the Taylor coefficients, here, in order to obtain the bounds required for the weighted summability, we need restrictions on the growth of the derivatives of the diffusion coefficient and right-hand side with respect to the parameter. 
\begin{theorem}[Weighted $\ell^2$-summability]\label{thm:l2summability} 
Let Assumption \ref{ass:afanalytic} hold and let there exist a positive sequence $(\rho_j)_{j\geq 1}$ and $0<\Kt<1$ such that
	\begin{equation}
		\sup_{\x\in D}\sum_{j \geq 1} \rho_j b_j\left(\x \right) \leq \Kb,  \label{eq:Kbassumption}
	\end{equation}
	 where the function sequence $(b_j(\x))_{j \geq 1}$ is such that, for a.e. $\x\in D$, $b_j(\x)\geq 0$  and the derivatives of $A$ and $F$ are bounded by 
	\begin{align}
		|\partial^\mu A(\bm{0},\x)|_{2,2} \leq ((b_j(\x))_{j \geq 1})^\mu f_A(|\mu|), \qquad \left|\partial^\mu F(\bm{0},\x)\right| \leq ((b_j(\x))_{j \geq 1})^\mu f_F(|\mu|),\label{eq:AFder}
	\end{align}		
	for all $|\mu| > 0$. Here, $f_A,f_F:\mathbb{N} \to \mathbb{R}^+$ are such that the series
	\begin{align}
		g_A(x):=\sum_{n \geq 1} \frac{f_A(n)}{n!}x^n \label{eq:gAdef}  
	\end{align} 
	has a convergence radius $\rho_A$, and $\frac{f_F(n)}{n!}$ has a monotonic, non decreasing majorant $f^*_F(n)$, such that the series
	\begin{align}
		g_F(x):=\sum_{n \geq 1} f^*_F(n)^2x^n \label{eq:gFdef}  
	\end{align} 
	has a convergence radius $\rho_F$. Moreover, let $\Kt$ be such that 
	\begin{equation}
		\Kt < \min\{\rho_A\rho_F, \rho_A\}, \label{eq:rhoArhoFass}
	\end{equation}
	and
	\begin{equation}
		g_A(\Kt)  < A_{min}. \label{eq:assKbKc}
	\end{equation}
	Then, we have
	\begin{equation}
		\sum_{\mu \in \mathcal{F}} \rho^{2\mu} \anorm{t_\mu}^2 \leq B < \infty, \label{eq:weightedsumresult}
	\end{equation}
	for $B=B(\Kt, g_A, g_F, A_{min}, D)$. 
\end{theorem}

\begin{proof}
This proof follows the general ideas put forward in \cite{Bachmayr2017}, adapted with different intermediate bounds. 
To aid our efforts, we define $\sigma_k$ as the weighted sum of all $k$-th order Taylor coefficients with positive weights $(\rho_j)_{j\geq 1}$, $k \geq 0$:
\begin{align}
	\sigma_k := \sum_{\mu \in \Lambda_k} \rho^{2\mu} \anorm{t_\mu}^2,
\end{align}
where $\Lambda_k:=\{\mu\in\mathcal{F}, \, |\mu| = k\}$.  

For $\sigma_k$, $k \geq 0$, we have, by Lemma \ref{col:taylorvar},
\begin{align}
	\sigma_k &\leq \underbrace{ A_{min}^{-1} \int_{{D}} \sum_{\mu \in \Lambda_k} \sum_{\nu \in S_\mu}     \frac{\rho^{2\mu}}{(\mu - \nu)!} \left|\A \nabla  t_\nu \right|   \abs{\nabla t_\mu} \dx}_{:=(\text{I})}   + \underbrace{A_{min}^{-1} \int_{{D}} \sum_{\mu \in \Lambda_k}  \frac{\rho^{2\mu}}{\mu!} |\partial^\mu {F}(\bm{0})||t_\mu| \dx}_{:=(\text{II})}.   \nonumber 
\end{align}
First, we bound $(\text{I})$ with the help of the Cauchy-Schwarz inequality:
\begin{align}
	(\text{I})
	&\leq A_{min}^{-1} \int_{{D}} \sum_{\mu \in \Lambda_k}\sum_{\nu \in S_\mu}     \frac{\rho^{\mu - \nu}}{(\mu - \nu)!} |\A|_{2,2} \left( \rho^\nu\abs{\nabla  t_\nu} \right) \left( \rho^\mu \abs{\nabla t_\mu}\right) \dx \nonumber\\
	&\leq A_{min}^{-1} \int_{{D}} \sum_{\mu \in \Lambda_k} \left(\sum_{\nu \in S_\mu}   \varepsilon(\mu, \nu)  \left( \rho^\nu\abs{\nabla  t_\nu}\right) ^2\right)^{\frac{1}{2 }} \left(\sum_{\nu \in S_\mu}    \varepsilon(\mu, \nu)  \left(\rho^\mu \abs{\nabla t_\mu}\right)^2\right)^{\frac{1}{2 }} \dx, \label{eq:firstCS}
\end{align}
with 
\begin{equation*}
	\varepsilon(\mu, \nu) := \frac{\rho^{\mu - \nu}|\A|_{2,2}}{ (\mu - \nu)!}.
\end{equation*}
To bound $\sum_{\nu \in S_\mu} \varepsilon(\mu, \nu)$, we first define $S_{\mu, l} := \{\nu \in S_\mu\,:\, \abs{\nu - \mu} = l \}$ for $l\geq 1$ and, using \eqref{eq:AFder} (and omitting explicit $\x$-dependence), we estimate:
\begin{align}
	\sum_{\nu \in S_{\mu, l}} \varepsilon(\mu, \nu) &= \sum_{\nu \in S_{\mu, l}} \frac{\rho^{\mu - \nu}|\A|_{2,2}}{ (\mu - \nu)!}
	\leq f_A(l) \sum_{\nu \in S_{\mu, l}} \frac{\rho^{\mu - \nu} ((b_j)_{j \geq 1})^{\mu - \nu} }{ (\mu - \nu)!} \nonumber\\
	&=  \frac{f_A(l)}{l!} \left( \sum_{j \geq 1} \rho_j b_j     \right)^l   \leq   \frac{f_A(l)}{l!} \Kt^l,\label{eq:enddouble}
\end{align}
where, in the last step, we have used equation \eqref{eq:Kbassumption} in the assumptions. From this, remembering that $\mu\in\Lambda_k$, we can bound the complete sum as
\begin{align}
	\sum_{\nu \in S_{\mu}} \varepsilon(\mu, \nu) &= \sum_{l=1}^k \sum_{\nu \in S_{\mu, l}} \varepsilon(\mu, \nu) 
	\leq  \sum_{l=1}^\infty \frac{f_A(l)}{l!} \Kb^l 
	= g_A(\Kt)
	= \deltatilde A_{min},   \label{eq:epsilonsumineq}
\end{align}
where $\deltatilde :=  \frac{g_A(\Kt)}{A_{min}} < 1$ by \eqref{eq:assKbKc}.
We insert estimate \eqref{eq:epsilonsumineq} into equation \eqref{eq:firstCS} and apply the Cauchy-Schwarz inequality to obtain
\begin{align}
	(\text{I}) &\leq  \sqrt{ \deltatilde  A_{min}^{-1}} \int_{{D}} \left(\sum_{\mu \in \Lambda_k}\sum_{\nu \in S_\mu}    \varepsilon(\mu, \nu)  \left( \rho^\nu\abs{\nabla  t_\nu}\right)^2 \right)^{\frac{1}{2 }}\left( \sum_{\mu \in \Lambda_k} \left(\rho^\mu \abs{\nabla t_\mu}\right)^2\right)^{\frac{1}{2 }} \dx. \label{eq:secondCS}
\end{align}
In order to treat the first double sum, we introduce, with $\nu \in \Lambda_l$ for $l\leq k-1$ and $k \geq 1$, $R_{\nu, k} = \{\mu \in \Lambda_k\, :\, \nu \in S_{\mu}\}$, such that, analogously to the estimates for $S_{\mu, l}$,
\begin{align}
	\sum_{\mu \in R_{\nu, k}} \varepsilon(\mu, \nu) 
	\leq \frac{f_A(k-l)}{(k-l)!}\Kb^{k-l}.\label{eq:secondepsineq}
\end{align}
Using inequality \eqref{eq:secondepsineq} and equation \eqref{eq:secondCS} leads to
\begin{align}
	(\text{I})  &\leq \sqrt{\deltatilde  A_{min}^{-1}}	\int_{{D}} \left(\sum_{l=0}^{k-1}\frac{f_A(k-l)}{(k-l)!}\Kb^{k-l} \sum_{\nu \in \Lambda_l}   \left(\rho^\nu\abs{\nabla  t_\nu}\right)^2\right)^{\frac{1}{2}}  \left(\sum_{\mu \in \Lambda_k}\left( \rho^\mu \abs{\nabla t_\mu}\right)^2\right)^{\frac{1}{2 }} \dx \nonumber\\
	&\leq \sqrt{\deltatilde  A_{min}^{-1}} \left(\sum_{l=0}^{k-1}\frac{f_A(k-l)}{(k-l)!}\Kb^{k-l} \sum_{\nu \in \Lambda_l} \int_{{D}} \left( \rho^\nu\abs{\nabla  t_\nu}\right)^2\dx \right)^{\frac{1}{2 }} \left(  \sum_{\mu \in \Lambda_k}\int_{{D}} \left( \rho^\mu \abs{\nabla t_\mu}\right)^2\dx\right)^{\frac{1}{2}} \nonumber\\
	&=  \sqrt{\deltatilde  A_{min}^{-1}} \left(\sum_{l=0}^{k-1}\frac{f_A(k-l)}{(k-l)!}\Kb^{k-l} \sigma_l \right)^{\frac{1}{2 }}  \sigma_k^{\frac{1}{2}}.\label{eq:Ibound}
\end{align}
The term $(\text{II})$ can by bounded via the Cauchy-Schwarz inequality and using \eqref{eq:AFder} (and that $\lVert\,\cdot\,\rVert_{\ell^2}\leq\lVert\,\cdot\,\rVert_{\ell^1}$), obtaining
\begin{align}
	 (\text{II})
	 & \leq A_{min}^{-1}\int_{{D}} \left(\sum_{\mu \in \Lambda_k} \left[ \frac{\rho^{\mu}}{\mu!} |\partial^\mu {F}(\bm{0})|\right]\right)\left(\sum_{\mu \in \Lambda_k}\rho^{2\mu}|t_\mu|^2 \right)^{\frac{1}{2}} \dx\nonumber\\
	 &\leq   A_{min}^{-1}f^*_F(k) \Kb ^k  \int_{{D}}    \left(\sum_{\mu \in \Lambda_k}\rho^{2\mu}|t_\mu|^2 \right)^{\frac{1}{2}} \dx.  \nonumber
\end{align}
Finally, to complete the bound of (II), we apply the Cauchy-Schwarz and Poincar\'e inequalities, and we use the assumption on $f_F$, resulting in
\begin{align}
	  (\text{II}) &\leq  A_{min}^{-1}f^*_F(k) \Kb^k |D|^{\frac{1}{2}} \Cpc  \left(\sum_{\mu \in \Lambda_k}\rho^{2\mu}\int_{{D}}|\nabla t_\mu|^2 \dx \right)^{\frac{1}{2}}
	  \leq A_{min}^{-1} f^*_F(k) \Kb^k   |D|^{\frac{1}{2}}   \Cpc   \sigma_k^{\frac{1}{2}},  \label{eq:IIbound}
\end{align}
where $C_p$ denotes the Poincar\'e constant, depending on the domain $D$. Now, combining equations \eqref{eq:Ibound} and \eqref{eq:IIbound} yields
\begin{align}
	\sigma_k &\leq    \Bigg(\sqrt{ \deltatilde  A_{min}^{-1}}  \left( \sum_{l=0}^{k-1}\frac{f_A(k-l)}{(k-l)!}\Kb ^{k-l}  \sigma_l\right)^{\frac{1}{2}} +  A_{min}^{-1} f^*_F(k) \Kb^k   |D|^{\frac{1}{2}}   \Cpc     \Bigg)^2.  \label{eq:squaredsigmakbound}
\end{align}
To show the final part of the proof, we will bound $\sigma_k$, $k\geq 0$, by 
\begin{equation}
	\sigma_k \leq \Csigma (k+1) \delta^k f^*_F(k)^2,  \label{eq:sigmaineq}
\end{equation} 
where $\Csigma := \Upsilon  \max(\sigma_0, A_{min}^{-2} |D|\Cpc^2)$ with  $\Upsilon\in\mathbb{R}$ such that
\begin{align}
	 \Upsilon \geq \max\left\{\frac{4\deltatilde }{\left(1-\deltatilde\right)^2}, 1\right\}, \label{eq:Nepsbound}
\end{align}
and where we choose $\delta$ such that $\frac{\Kt}{\rho_A} < \delta < \min\{\rho_F,1\}$,
the existence of such $\delta$ being ensured by the hypothesis \eqref{eq:rhoArhoFass}.

We proceed by induction. The base case for $k=0$ is shown by noting that $\sigma_0 \leq \Csigma$.
The inductive step is proven by expansion of equation \eqref{eq:squaredsigmakbound}. For $k\geq 1$:
\begin{align*}
	\sigma_k &\leq       \Bigg( \sqrt{ \deltatilde A_{min}^{-1}}  \left( \Csigma f^*_F(k)^2\delta^k k  \sum_{l=0}^{k-1}\frac{f_A(k-l)}{(k-l)!} \left(\frac{\Kb }{\delta}\right)^{k-l} \right)^{\frac{1}{2}} + A_{min}^{-1} f^*_F(k) \Kb^k   |D|^{\frac{1}{2}}   \Cpc   \Bigg)^2\nonumber\\
	&\leq    \left(    \left( \deltatilde\Csigma f^*_F(k)^2\delta^k k\right)^\frac{1}{2} +  A_{min}^{-1} f^*_F(k) \delta^k   |D|^{\frac{1}{2}}   \Cpc   \right)^2\\
	&\leq  \deltatilde  \Csigma f^*_F(k)^2\delta^k k +A_{min}^{-2} f^*_F(k)^2 \delta^{k} |D|\Cpc^2 + 2 \deltatilde^{\frac{1}{2}} \Csigma^\frac{1}{2}f^*_F(k)^2   \delta^k k A_{min}^{-1}   |D|^\frac{1}{2}\Cpc,
\end{align*}
where we remind that the last inequality follows from $\delta<1$. Now, we simplify the term $\Csigma^\frac{1}{2} A_{min}^{-1}   |D|^\frac{1}{2}\Cpc$ by bounding it by $\frac{\Csigma}{\sqrt{\Upsilon}}$, and using that $\Upsilon \geq 1$:
\begin{align*}
	\sigma_k &\leq  \deltatilde  \Csigma f^*_F(k)^2\delta^k k + \frac{\Csigma f^*_F(k)^2}{\Upsilon} \delta^{k} + \frac{2\deltatilde^{\frac{1}{2}} }{\sqrt{\Upsilon}} \Csigma f^*_F(k)^2\delta^k k \\
	&\leq  \Csigma \delta^k f^*_F(k)^2\left( \left(\deltatilde + \frac{2\deltatilde^{\frac{1}{2}} }{\sqrt{\Upsilon}} \right)k+1\right), 
\end{align*}
which advances the induction due to \eqref{eq:Nepsbound}.
Finally, we are able to conclude our proof by summing over $k$ and employing equation \eqref{eq:sigmaineq}:
\begin{equation*}
	\sum_{k\geq 0} \sigma_k \leq \Csigma  \sum_{k\geq0}(k+1) f^*_F(k)^2\delta^k =\Csigma \left(g_F(\delta)+\delta g_F'(\delta) \right)< \infty,
\end{equation*}
which converges by the choice of $\delta$.
\end{proof}
\begin{remark}\label{rem:seperation1}
	As we will see also further in Section \ref{sec:modelproblem}, what matters is not much the maximal amount of variation (with respect to the parameter) in the coefficient $A$ per se, encoded in the constant $\Kb$, but rather the amount of variation relative to $A_{min}$ and its relationship with the amount of variation in the right-hand side, encoded in the relationships \eqref{eq:assKbKc} and \eqref{eq:rhoArhoFass}, respectively. (although the bound $B$ in \eqref{eq:weightedsumresult} can depend on $\Kb$, and the relationship between $A$ and $F$). 
\end{remark}
\begin{remark}
	Theorem \ref{thm:l2summability} introduces a balance between $\rho_A$ and $\rho_F$, where, for given $\Kt$, large $\rho_A$ allows for small $\rho_F$ and vice-versa. This inverse relationship is not surprising, as it follows from the underlying PDE structure.
\end{remark}

\subsection{Summability of Taylor coefficients}
As a direct consequence of the weighted summability proved in Theorem \ref{thm:l2summability}, we have the following $\ell^p$-summability result:
\label{sec:ellpsum}
\begin{theorem}[$\ell^p$-summability]
	\label{thm:lpsummability}
	Let the assumptions of Theorem \ref{thm:l2summability} hold with $\rho_j > 1$, $j\geq 1$, and $(\rho_j^{-1})_{j\geq 1} \, \in  \, \ell^q(\mathbb{N})$, where $q = \frac{2p}{2-p}$ for some $p < 2$. Then we have $(\anorm{t_\mu})_{\mu \in \mathcal{F}} \, \in \, \ell^p(\mathcal{F})$. 
\end{theorem}
\begin{proof}
	The proof follows closely the one of Corollary 2.3 in \cite{Bachmayr2017a}.	By Theorem \ref{thm:l2summability} and H\"older's inequality,
	\begin{equation}
		\sum_{\mu \in \mathcal{F}} \anorm{t_\mu}^p \leq \left( \sum_{\mu \in \mathcal{F}} \rho^{2\mu} \anorm{t_\mu}^2 \right)^{\frac{p}{2}} \left( \sum_{\mu \in \mathcal{F}} \rho^{-\frac{2p}{2-p}\mu} \right)^{\frac{2-p}{2}}\leq B \left( \sum_{\mu \in \mathcal{F}} \rho^{-\frac{2p}{2-p}\nu} \right)^{\frac{2-p}{2}},
	\end{equation}		
	with $B$ as in \eqref{eq:weightedsumresult}. Moreover,
	\begin{equation}
		\sum_{\mu \in \mathcal{F}} \rho^{-\frac{2p}{2-p}\mu} = \prod_{j \geq 1}\left( \sum_{k=0}^\infty \rho_j^{-qk}  \right)= \prod_{j \geq 1}\left( 1-\rho_j^{-q} \right)^{-1},
	\end{equation}
	and the last product converges precisely when $(\rho_j^{-1})_{j\geq 1} \, \in  \, \ell^q(\mathbb{N})$.
\end{proof}
From the $\ell^p$-summability of the Taylor coefficients of the solution, convergence rates for truncated Taylor expansions follow:
\begin{corollary}[Best $N$-term approximation]\label{cor:stechkin}
	Let the assumptions of Theorem \ref{thm:lpsummability} hold. Then, for the best $N$-term approximation
	\begin{align*}
		u_N^T(\y) := \sum_{\mu \in \Lambda_N^T} t_\mu \y^\mu,
	\end{align*}
	where $\Lambda_N^T$ is the index set corresponding to the $N$ Taylor coefficients with largest $V$-norm, it holds that
	\begin{align*}
		\|u-u_N^T\|_{L^\infty(Y,V)} \leq C (N+1)^s, \,\,\,\, s:=\frac{1}{p} - 1 = \frac{1}{q} - \frac{1}{2}, \quad C:= \|(\anorm{t_\mu})_{\mu \in \mathcal{F}}\|_{\ell^p(\mathcal{F})}.
	\end{align*}
\end{corollary}
\begin{proof}
This results, as usual, from the application of Stechkin's lemma \cite{devore1998,Gogoladze2022}, combined with the summability properties of Theorem \ref{thm:lpsummability}.
\end{proof}
\begin{remark}[Comparison with \cite{Bachmayr2017a}]
	Theorem \ref{thm:lpsummability} is a generalization of Theorem 1.2 in \cite{Bachmayr2017a} to non-affine parameter dependence (and parameter-dependent right-hand side). Specifically, when taking $A(\y)= a(\y)I_d = \left(\bar{a} + \sum_{j \geq 1} y_j \psi_j\right)I_d$, $I_d$ being the $d\times d$ identity matrix, and $F$ independent of $\y$, the main result in \cite{Bachmayr2017a} is recovered by choosing $b_j=|\psi_j|$, $f_A(n)=1$ for $n=1$ and $f_A(n)=0$ otherwise, and $f_F(n)\equiv 0$. The difference between the two results is that, in the present work, $g_A(x)$ in \eqref{eq:gAdef} is more general than for the affine case, where $g_A(x)=x$. For this reason, \eqref{eq:assKbKc} gives a slightly stronger condition than the one needed in \cite{Bachmayr2017a}, referred there as \emph{(UEA\textsuperscript{$\ast$})} assumption.
\end{remark}

\medskip
\begin{remark}[Sharpness of Theorem \ref{thm:lpsummability}]
	On the basis of the previous remark, the sharpness of our summability result in Theorem \ref{thm:lpsummability} follows from the example in \cite[Sect. 4.2]{Bachmayr2017a}.
\end{remark}

%% file: 4application/4application.tex
\section{Application to parameterized domains}
\label{sec:modelproblem}
To show the relevance of Theorem \ref{thm:lpsummability}, we consider the elliptic diffusion equation on family of open Lipschitz domains $\mathcal{D}(\y) \subset \mathbb{R}^d$, parameterized by $\bm{y} \in Y$: 
\begin{equation}
	\begin{cases}
		-\nabla\cdot \left(a\nabla \uparamdomain \right) = f & \text{in } \mathcal{D}(\y),\\
		\uparamdomain = 0 & \text{on } \partial \mathcal{D}(\y),\\
		\text{for every  }\bm{y} \in Y,
	\end{cases} \label{eq:shapepde}
\end{equation}
where $a\in L^\infty(\mathcal{D}_H), f\in L^2(\mathcal{D}_H)$ are analytic functions from the hold-all domain $\mathcal{D}_H = \cup_{\y\in Y}\mathcal{D}(\y)$ to $\mathbb{R}$, fixed and independent of $\y\in Y$. We take care to distinguish $\uparamdomain$, the solution on the parameterized domain, from $u$, which later will be the solution on a fixed domain. 
Finally, we assume $a$ satisfies the following regularity assumption:
\begin{assumpttion}\label{ass:a}
	Both $a$ and $f$ in \eqref{eq:shapepde} are analytic functions from $\mathcal{D}_H$ to $\mathbb{R}$. Next to this, we have that the coefficient $a$ is uniformly bounded by 
	\begin{equation}
		0 < a_{min} \leq a(\x) \leq a_{max} < \infty,
	\end{equation}
	for some constants $a_{min},a_{max}>0$.
\end{assumpttion}

For each $\y \, \in \, Y$, the variational formulation of \eqref{eq:shapepde} reads:
\begin{align}
&\text{find }\uparamdomain \in H_0^1(\Dy)  \text{ such that}\nonumber \\
	 &\qquad\int_{\Dy} a \nabla \uparamdomain \cdot \nabla v \d \bm{x} = \int_{\Dy} f v\d \bm{x}, \quad \text{ for all } {v} \in V. \label{eq:variational_formulation_y}
\end{align}
Thanks to the assumptions on $a$ and $f$, well-posedness for every $\y \in Y$ is ensured by the Lax-Milgram lemma. 

We model the parameterized domain as a deformation of a Lipschitz reference domain, which we call ${D}$. In the spirit of \cite{Castrillon-Candas2016,Harbrecht2016,Hiptmair2018}, we apply a mapping approach \cite{Tartakovsky2006,Xiu2006} to pull the variational formulation \eqref{eq:variational_formulation_y} defined on $\mathcal{D}(\y)$ back onto the reference domain ${D}$ using a parameter-dependent diffeomorphism $\Phi(\y) :{D} \to \mathcal{D}(\y)$.
Following \cite{Cohen2018}, and motivated by similar approaches in shape optimization \cite{Haubner2020}, we consider maps that depend affinely on the entries of $\y$:
\begin{equation}
	\Phi(\y; \x) = \x + \sum_{j\geq 1} \Phi_j({\bm{x}})y_j,\quad \x \in {D}.\label{eq:phidef}
\end{equation}
For these maps, we formulate the following assumption:
\begin{assumpttion}\label{ass:phij}
	The sequence $(\Phi_j)_{j\geq 1}$ satisfies the following smoothness and decay properties:
	\begin{itemize}
		\item $\Phi(\y;\cdot)$,  $\Phi_j(\cdot)$  $\in W^{1, \infty}(D)$,  for all $j \geq 1$, and $\Phi^{-1}(\y;\cdot) \in W^{1, \infty}(\mathcal{D}(\y))$ for all $\y \in Y$;
		\item there exist constants $C_\Phi, C_{\Phi^{-1}}$, independent of $\y \in Y$, such that we have the bounds \\$\|\Phi(\y;\cdot)\|_{ W^{1, \infty}(D)}$ $< C_\Phi$ and $\|\Phi^{-1}(\y;\cdot)\|_{W^{1, \infty}(\mathcal{D}(\y))} < C_{\Phi^{-1}}$ for all $\y \in Y$.
	\end{itemize}
\end{assumpttion}
We can take derivatives of $\Phi(\y,\x)$ to obtain the Jacobian matrix $\D\Phi(\y;\x)$. 
As a consequence of the the Courant-Fisher theorem for singular values \cite[Thm. 3.1.2]{Horn1991}, applied to $\D\Phi(\y)$ and $\D\Phi^{-1}(\y)$, we have:
\begin{lemma} \label{lem:courantfish}
	Let Assumption \ref{ass:phij} hold. We have that, for the singular values of $\D\Phi^{-1}(\y)$, $\sigma_1=\sigma_1(\y; \x), \cdots, \sigma_d=\sigma_d(\y; \x)$, there exist constants $\sigma_{min},  \sigma_{max} > 0$, independent of $\y \in Y$, such that 
	\begin{equation}
		\sigma_{min} \leq \sigma_1(\y; \x), \cdots, \sigma_d(\y; \x) \leq \sigma_{max}, \quad  \text{for a.e. }\x\in \mathcal{D}(\y)\text{ and for all }\y \in Y .
	\end{equation}
\end{lemma}
With the mapping $\Phi(\y; \x)$ at hand, we can write the following variational problem on the reference domain:
\begin{align}
	&\text{for every }\y \in Y,\text{ find }{u}(\y) \in V \text{ such that} , \nonumber\\
	&\qquad \int_{{D}} {A}(\y) {\nabla} {u}(\y) \cdot  {\nabla} {v} \d\x = \int_{{D}} F(\y) {v} \dx, \quad \text{for all } {v} \in V, \label{eq:transformedvarform}
\end{align}
where $V=H_0^1(D)$,
\begin{align}
	 {A}(\y; \x) := (a\circ\Phi(\y; \x)) \D\Phi^{-1}(\y; \x) \D\Phi^{-\top}(\y; \x)\det(\D\Phi(\y; \x)), \label{eq:hatAdef}
\end{align}
and 
\begin{align}
	 F(\y; \x) := (f\circ\Phi(\y; \x)) \det(\D\Phi(\y; \x)).\label{eq:hatfdef}
\end{align}
Under Assumptions \ref{ass:a} and\ref{ass:phij}, the problem \eqref{eq:transformedvarform} is well-posed and $\|u(\y)\|_{V}$ is bounded independent of $\y \in Y$.

To efficiently compute evaluations of ${u}(\y)$, we are interested in a polynomial approximation for the mapping $\y \mapsto {u}(\y)$ from $Y$ to $V$. The holomorphy of an extension of the map $\y \mapsto {u}(\y)$  from $Y$ to $V$ to complex poly-ellipses, under Assumptions \ref{ass:a} and \ref{ass:phij}, was shown in \cite{Chkifa2015,Cohen2018,Hiptmair2018}. Therefore, we can expand the solution to \eqref{eq:transformedvarform} in terms of a Taylor series \eqref{eq:tayloreq}. We will investigate the summability properties of the Taylor coefficients by employing Theorem \ref{thm:lpsummability} in the following section.

\subsection{Setup for Theorem \ref{thm:lpsummability}}
\label{sec:verificationthm}
By Assumptions \ref{ass:a}, \ref{ass:phij}, and Lemma \ref{lem:courantfish}, Assumption \ref{ass:afanalytic} is satisfied. To apply Theorem \ref{thm:lpsummability}, we need a sequence $(b_j(\x))_{j \geq 1}$ that verifies its assumptions. In this subsection, we will first show   \textsl{pointwise} bounds on $|{A}(\bm{0}; \x)|_{2,2}$ and $|F(\bm{0}; \x)|$ from \eqref{eq:hatAdef}--\eqref{eq:hatfdef} in terms of  $ |\D \Phi_j(\bm{0}; \x) |_{2,2}$ and $ | \Phi_j (\bm{0}; \x)|$, $j\geq 1$, and use these bounds to choose $(b_j(\x))_{j \geq 1}$ appropriately. 
We define the following bounds on the spectral norms of the total derivatives of $a$ and $f$:
\begin{align*}
	\Ca(n):=\sup_{\x \in \mathcal{D}_H}|\D^na|_{n,2}, \quad \Cf(n):=\sup_{\x \in \mathcal{D}_H}|\D^n f|_{n,2}.
\end{align*}
Moreover, we let $\Cat$ and $\Cft$ be the monotonically non-decreasing majorants of $\Ca$ and $\Cf$ respectively.

Now, the desired bounds on $|A(\bm{0}; \x)|_{2,2}$ and $|F(\bm{0}; \x)|$ follow from direct computation and are given by the following lemma:
\begin{lemma}
	\label{thm:Afbound}
	For $\y \in Y $, let ${A}(\y)$ and $F(\y)$ be defined by equations \eqref{eq:hatAdef} and \eqref{eq:hatfdef}, respectively, and let Assumption \ref{ass:a} and \ref{ass:phij} hold. Then, for a.e. $\x \in {D}$ and for every $\mu \in \mathcal{F}$, we have the bounds
	\begin{align*}
		&|\partial^\mu{A}(\bm{0}; \x)|_{2,2} \leq \\
		&\quad\Cat(|\mu|)  \sum_{ \nu \leq \mu} \binom{\mu}{\nu} \frac{(|\mu| - |\nu| + 2)!}{2} ((|\D \Phi_j(\bm{0}; \x)|_{2,2} )_{j \geq 1})^{\mu - \nu} \sqrt{2}^{\,|\mu - \nu|} ((| \Phi_j(\bm{0}; \x)| )_{j \geq 1})^{ \nu},
	\end{align*}
	and 
	\begin{equation*}
		|\partial^\mu F(\bm{0}; \x)| \leq \Cft(|\mu|)   \sum_{\nu \leq \mu} \binom{\mu}{ \nu } |\nu |!\sqrt{2}^{|\nu |}((|\D \Phi_j(\bm{0}; \x)|_{2,2})_{j\geq 1})^\nu   ((|\Phi_j(\bm{0}; \x)|)_{j \geq 1})^{\mu - \nu}.
	\end{equation*}
\end{lemma}
\begin{proof}
We expand the derivative  $\partial^{\mu}{A}(\bm{y})$ using the general Leibniz rule and substitute $\y = \bm{0}$ to obtain:
\begin{align}
	&|\partial^\mu {A}(\bm{0}) |_{2,2} \leq \sum_{\nui+\nuii+\nuiii+\nuiv=\mu} \binom{\mu}{\nu^{(1)},\nu^{(2)},\nu^{(3)}, \nuiv}\,\,\,\, \left| \partial^{\nu^{(1)}}\left[ \left( \D \Phi \right)^{-1} (\bm{0})\right]\right|_{2,2} \nonumber\\
&\qquad \qquad\cdot\left|   \partial^{\nu^{(2)}}\left[ \left( \D \Phi \right)^{-\top} (\bm{0})\right]   \right|_{2,2}\cdot\left|\partial^{\nu^{(3)}}\left[\det\left( \D \Phi \right) (\bm{0})\right]\right|_{2,2}\cdot | \partial^{\nuiv}\left[\left( a \circ \Phi \right)(\bm{0}) \right]|.\label{eq:inbetweenbound}
\end{align}
This leaves us with the task to estimate each norm on the right-hand side of the previous equation. First, we tackle $\left| \partial^\mu \left[ \left( \D \Phi \right)^{-1}(\bm{0})\right]\right|_{2,2}$, by computing the derivative explicitly. We do so instead of employing the Fa\`a di Bruno formula to exploit the affine expansion. We obtain:
\begin{align*}
	\partial^\mu \left[ \left( \D \Phi \right)^{-1}(\bm{y})\right] &= (\D \Phi)^{-1}(\y) \sum_{\xi\in P(\mu)} \prod_{j=1}^{|\mu|}\left[ \D \Phi_{\xi_j} (\D \Phi)^{-1}(\y) \right], 
\end{align*}
where $P(\mu)$ is the set of all possible permutations of the derivatives in $\mu$, each counted with its multiplicity. We can evaluate this at $\bm{y}=\bm{0}$, where $\D\Phi(\bm{0})=I_d$,  to obtain:
\begin{align}
	\left| \partial^\mu \left[ \left( \D \Phi \right)^{-1}(\bm{0})\right]\right|_{2,2} &\leq  \sum_{\xi\in P(\mu)} \prod_{j=1}^{|\mu|} |\D \Phi_{\xi_j} |_{2,2}
	=  \abs{\mu}! ((|\D \Phi_j |_{2,2})_{j \geq 1}) ^{\mu}.\label{eq:invbound}
\end{align}
Similarly, for the transpose, we obtain:
\begin{align}
	\left| \partial^\mu \left[ \left( \D \Phi \right)^{-\top}(\bm{0})\right]\right|_{2,2} &\leq  \abs{\mu}! \left((|\D \Phi_j |_{2,2})_{j \geq 1}\right) ^{\mu},\label{eq:invtranspbound}
\end{align}
because the spectral norm is transpose invariant.
Next, we take a look at the determinant by evaluating the derivative iteratively. For an arbitrary $j \in \text{supp}(\mu)$:
\begin{align*}
	g_\mu^{\bm{y}}&:=\left| \partial^\mu \left[ \det\left( \D \Phi(\y) \right)\right]\right| \\
	&\phantom{:}=\left| \partial^{\mu - e_j}\left[ \det(\D \Phi(\y)) \tr\left((\D\Phi(\y)) ^{-1}\D\Phi_j\right)\right]\right|\\
	&\phantom{:}=\left| \sum_{\nu \leq \mu - e_j}  \binom{\mu - e_j }{ \nu}  \partial^{\nu} \left[ \det\left( \D \Phi(\y) \right)\right] \tr\left( \partial^{\mu - e_j - \nu}\left[ (\D \Phi(\y))^{-1} \right] \D \Phi_j  \right)\right|\\
	&\phantom{:}\leq \sum_{\nu \leq \mu - e_j}  \binom{\mu - e_j }{\nu}  g_\nu^{\bm{y}} | \partial^{\mu - e_j - \nu}\left[ (\D \Phi(\y))^{-1} \right]|_F | \D \Phi_j |_F,
\end{align*}
with $|\cdot|_F$ the Frobenius norm.
Evaluating at $\bm{y}=\bm{0}$ and employing \eqref{eq:invbound}, we obtain
\begin{align*}
	\left| \partial^\mu \left[ \det\left( \D \Phi(\bm{0}) \right)\right]\right| = g_\mu^{\bm{0}}&\leq \sum_{\nu \leq \mu - e_j}  \binom{\mu - e_j }{ \nu}  g_\nu^{\bm{0}} (|\mu - \nu|-1)!((| \D \Phi_j |_F)_{j\geq 1})^{\mu - \nu},
\end{align*}
and we claim this recurrence relation is solved by 
\begin{equation}
	g_\mu^{\bm{0}} \leq \abs{\mu}!(|\D \Phi_j|_F)^\mu. \label{eq:detbound}
\end{equation}
To show this, we proceed by induction on $\abs{\mu}$. For $\mu = \bm{0}$, we have  $g_{\bm{0}}^{\bm{0}} = 1=|\bm{0}|!(|\D \Phi_j|_F)^{\bm{0}},$ where we used the convention that $0!=1$.
Assuming \eqref{eq:detbound} holds for $\abs{\mu} \leq n$, for $g_{\mu+e_j}^{\bm{0}}$, and any $j\geq 1$, we have
\begin{align*}
	g_{\mu+e_j}^{\bm{0}}
	&\leq((| \D \Phi_j |_F)_{j\geq 1})^{\mu}\sum_{\nu \leq \mu}  \binom{\mu }{ \nu}  \abs{\nu}! |\mu - \nu|!
	\\&=   (\abs{\mu}+1)!((| \D \Phi_j |_F)_{j\geq 1})^{\mu},
\end{align*}
where, in the last inequality, we used the case $N=2$ of Lemma \ref{lem:multinomialsum} in Appendix \ref{ap:combilemma}.

Finally, we bound $| \partial^\mu (a \circ \Phi) |$ by exploiting the affine dependence of $\Phi$ on $\y$:
\begin{align}
	| \partial^\mu ( a \circ \Phi ) | \leq | \D^{|\mu|} a|_{|\mu|,2}  ((|\Phi_j|)_{j \geq 1})^{\mu}\leq \Cat(|\mu|)((|\Phi_j|)_{j \geq 1})^{\mu}.\label{eq:gbound}
\end{align}

We finish the first part of the proof by expanding equation \eqref{eq:inbetweenbound} using the bounds in equations \eqref{eq:invbound} -- \eqref{eq:gbound} and we obtain:
\begin{align*}
	|\partial^\mu{A}(\bm{0})|_{2,2} 
	&\leq \sum_{\nu^{(1)}+\nu^{(2)}+\nu^{(3)}+\nu^{(4)}=\mu} \binom{\mu}{ \nu^{(1)},\nu^{(2)},\nu^{(3)}, \nuiv} |\nu^{(1)}|!\left( (|\D \Phi_j |_{2,2})_{j \geq 1} \right)^{\nu^{(1)}}\nonumber\\
& \qquad    |\nu^{(2)}|!\left( |\D \Phi_j |_{2,2})_{j \geq 1} \right)^{\nu^{(2)}}  |\nu^{(3)}|!\sqrt{2}^{|\nu^{(3)}|}\left( |\D \Phi_j |_{2,2})_{j \geq 1} \right)^{\nu^{(3)}} \Cat(|\nuiv|)((|\Phi_j|)_{j \geq 1})^{\nuiv} \nonumber \\
	\leq  &\sum_{\nuiv \leq \mu} \left[ \binom{\mu }{ \nuiv} \sqrt{2}^{\abs{\mu - \nuiv}} \left(( |\D \Phi_j |_{2,2})_{j\geq 1} \right)^{\mu - \nuiv}\right. \\ &\qquad \left. \sum_{\nu^{(1)}+\nu^{(2)}+\nu^{(3)}=\mu-\nuiv} \binom{\mu - \nuiv }{ \nu^{(1)},\nu^{(2)},\nu^{(3)}} |\nu^{(1)}|! |\nu^{(2)}|! |\nu^{(3)}|!  \Cat(|\nuiv|) ((|\Phi_j|)_{j \geq 1})^{\nuiv}\right]  \nonumber \\
	\leq  &\Cat(|\mu|)  \sum_{ \nuiv \leq \mu} \binom{\mu}{\nuiv} \frac{(|\mu| - |\nuiv| + 2)!}{2} ((|\D \Phi_j|_{2,2} )_{j \geq 1})^{\mu - \nuiv} \sqrt{2}^{\,|\mu - \nuiv|} ((| \Phi_j| )_{j \geq 1})^{ \nuiv},
\end{align*}
where, on the last line, we employed Lemma \ref{lem:multinomialsum} in Appendix \ref{ap:combilemma} with $N=3$.

To find the derivative $\left|\partial^\mu F(\bm{0}) \right|$, and finish the proof, we expand using the general Leibniz rule, such that equation \eqref{eq:detbound}, and the analogous for $f$ of \eqref{eq:gbound}, lead to
\begin{align*}
	\left|\partial^\mu F(\bm{0}) \right| &\leq \sum_{\nui + \nuii = \mu} \binom{\mu }{ \nui, \nuii} | \partial^\nui \det(\D\Phi(\bm{0})) |  | \partial^\nuii (f \circ \Phi(\bm{0})) |\\
	&\leq \Cft(|\mu|)   \sum_{\nu \leq\mu} \binom{\mu }{ \nu} |\nu|!\sqrt{2}^{|\nu|}(|(\D \Phi_j|_{2,2})_{j \geq 1})^\nu ((|\Phi_j|)_{j \geq 1})^{\mu-\nu}.
\end{align*}
\end{proof}
With these bounds at hand, we define $(b_j(\x))_{j \geq 1}$ in the following corollary:
\begin{corollary}
	\label{lm:Afbound}
	For $\y \in Y $, let ${A}(\y; \x)$, $F(\y; \x)$ be defined by equations \eqref{eq:hatAdef} and \eqref{eq:hatfdef} respectively, and let Assumptions \ref{ass:a} and \ref{ass:phij} hold. Then, for a.e. $\x \in {D}$ and $\mu \in \mathcal{F}$, we have the bounds
	\begin{align}
		|\partial^\mu A(\bm{0};\x)|_{2,2} \leq ((b_j(\x))_{j \geq 1})^\mu f_A(|\mu|), \qquad |\partial^\mu F(\bm{0};\x)| \leq ((b_j(\x))_{j \geq 1})^\mu f_F(|\mu|),
	\end{align}	
	with
	\begin{align}
		f_A(n):=\Cat(n) \frac{(n+3)!}{3!} \text{, } f_F(n):=\Cft(n) n!,
	\end{align}
	and
	\begin{equation}\label{eq:bj_phi}
		b_j(\x):= |\Phi_j (\x)| + \sqrt{2}|\D\Phi_j(\x) |_{2,2}.
	\end{equation}
\end{corollary}
\begin{proof}
	The result follows from the application of Lemma \ref{lem:multinomialsum} with $N=2$ to Lemma \ref{thm:Afbound}, after we bound both $\sqrt{2}|\D\Phi_j|_{2,2}$ and $|\Phi_j|$ by $b_j$, $j \geq 1$.
\end{proof}
To verify the assumptions of Theorem \ref{thm:lpsummability}, we are left with the task to verify equations \eqref{eq:rhoArhoFass} -- \eqref{eq:assKbKc} and to find a positive sequence $(\rho_j)_{j\geq 1}$ with $\rho_j > 1$ satisfying equation \eqref{eq:Kbassumption}. To attain this goal, we need to bound $b_j(\x)$, $j \geq 1$, in order to determine $\Kt < 1$ in \eqref{eq:Kbassumption}, to finally show 
\begin{equation}
	g_A(\Kb) = \sum_{n  \geq 1} \Cat(n) \binom{n+3}{3} \Kb^n 
	< a_{min}\sigma_{min}^4, \label{eq:gAnonexplicit}
\end{equation} 
such that, by Lemma \ref{lem:courantfish}, we have 
\begin{align}
	 g_A(\Kb) <  a_{min}\sigma_{min}^4  = A_{min} \text{ and } \max\{\rho_A\rho_F, \rho_A\} > \Kt.\label{eq:gAdefshape}
\end{align}
Here $\sigma_{min}$ is, following Lemma \ref{lem:courantfish}, the lower bound on the singular values of $\D\Phi^{-1}$.
To bound $b_j(\x)$, $j \geq 1$, we are dependent on bounding $\Phi_j(\x)$. To determine these partial transformations, we turn our attention to star-shaped domains in the next subsection.
\begin{remark}
	In the special case where $\Ca(n)$ is bounded, and therefore $\Cat(n)=\Cat$ is constant, we can expand equation \eqref{eq:gAnonexplicit} explicitly into 
	\begin{align*}
		g_A(\Kb) = \Cat \left( \frac{1}{(1-\Kt)^4} -1\right).
	\end{align*}
\end{remark}
\begin{remark}
	We observe that sharper constants in the convergence estimate could be obtained by adapting the proof of Theorem \ref{thm:l2summability} to the specific case of parametric domains considered in this section. This would improve the constants whilst leading to the same summability exponent for the Taylor coefficients.  
\end{remark}
\begin{remark}
	From equation \eqref{eq:gAdefshape}, we observe that, as described in Remark \ref{rem:seperation1}, the important quantity here is the amount of relative variation to $g_A(\Kt)/(a_{min}\sigma_{min}^4)$, where we see how different domain mappings have an effect on the constants.
\end{remark}

\subsection{Star-shaped domains}
\label{sec:starshapeddomains}
\label{sec:application}
Thus far, we have not considered an explicit expression for the parameterized domain $\Dy$.  We now illustrate how the results of the previous subsection apply to a star-shaped reference domain ${D}$ in $\mathbb{R}^2$. Due to this star-shape property, we express the latter in polar coordinates as the region inside a simple, closed curve $r_0(\theta)$, $\theta\in[0,2\pi)$:
\begin{align*}
	D = \{ \x \in \mathbb{R}^2: |\x| \leq r_0(\theta(\x)) \},
\end{align*} 
where $r_0 \in C^{0, 1}_{per}([0, 2\pi))$ parametrizes the boundary of the nominal domain ${D}$ and $\theta(\x)$ denotes the argument of $\x\in\mathbb{R}^2$.   
Similarly, we assume all parameterized domains $\mathcal{D}(\y)$ to be star-shaped as well.
Hence, we describe, for every $\y \in Y$, the parameterized domain $\Dy$ by 
\begin{align}
	\Dy= \{\x \in \mathbb{R}^2: |\x| \leq r(\y;\theta(\x))  \}. \label{eq:stardomain}
\end{align}
We model the boundary $r(\y;\theta)$ as an affine combination of the nominal radius and basis functions $\{\psi_j(\x)\}_{j \geq 1}$, namely
\begin{align}
	r(\y;\theta) = r_0(\theta) + \sum_{j \geq 1}\psi_j(\theta)y_j, \qquad\text{for all }\theta\in[0, 2\pi),\, \y\in Y. \label{eq:rdef}
\end{align}
To ensure, in the next subsections, the required regularity of the domain mapping $\Phi$ and its inverse, we add the following assumption and consecutive lemma:
\begin{assumpttion} \label{ass:phiregularity}
	For all $j \geq 1$, $\psi_j \in C^{0, 1}_{per}([0,\pi))$, and there exists $\gamma \in (0,1)$, independent of $\y\in Y$, such that
	\begin{equation}
		\sum_{j \geq 1}\psi_j(\theta)y_j \leq \gamma r_0(\theta), \qquad \sum_{j \geq 1}\psi_j'(\theta)y_j  < \infty, \qquad\text{for all }\theta\in[0, 2\pi).
	\end{equation}	
\end{assumpttion}
\begin{lemma}
	By Assumption \ref{ass:phiregularity}, we have
	\begin{equation}
		r(\y; \cdot) \in C^{0, 1}_{per}([0,2\pi)),\quad \text{for all } \y \in Y,
	\end{equation}
	and $\|r(\y, \cdot)\|_{C^{0,1}}$ has a  $\y$-independent upper bound.
\end{lemma}
With the parametrization \eqref{eq:rdef} at hand, we can find a bijective mapping $\Phi(\y;\cdot)$ from ${D}$ to $\Dy$, which depends linearly on $r(\y;\cdot)$. Due to this linear dependence and the affine expansion \eqref{eq:rdef}, we can express the domain transformation as in \eqref{eq:phidef}, and utilize it in Theorem \ref{thm:lpsummability}. 

One may use different techniques to obtain the mapping $\Phi$. One approach involves formulating an explicit radial transformation between ${D}$ and $\mathcal{D}(\y)$, by mapping the boundary of ${D}$ to the boundary of $\mathcal{D}(\y)$ using a mollifier. The use of such an explicit expression will offer several advantages when studying theoretical properties of this mapping. In Subsection \ref{sec:explicittransformation}, we will further explore this, referred to as mollifier mapping later.

Another approach is to describe the mapping implicitly as the solution of a linear partial differential equation, with \eqref{eq:rdef} as boundary condition.  
One such partial differential equation is Laplace's equation, applied to each component of the spatial coordinate separately. The resulting domain mapping is called the harmonic extension \cite{Li2001,Xiu2006}. The smoothing properties of the harmonic extension bring numerical advantages compared to the mollifier approach, as we will see in the numerical experiments. On the other hand, analytically proving the bounds required to apply Theorem \ref{thm:lpsummability} is more challenging, due to the lack of a general explicit expression of $\Phi$. In Subsection \ref{sec:harmonic}, we will delve further into the harmonic mapping.

The mollifier and the harmonic mappings, discussed in the upcoming subsections, allow us to illustrate, analogously to the findings in \cite{Bachmayr2017,Bachmayr2017a}, how the use of basis functions with localized supports can lead to faster convergence of the Taylor coefficients for the parametric solutions to \eqref{eq:pde}.

\subsection{Mollifier transformation}
\label{sec:explicittransformation}
In this subsection, we examine a transformation that maps the nominal domain ${D}$ to the parameterized domain $\Dy$ by employing a mollifier $\chi(\x)$. Following \cite{Hiptmair2018}, we define a mollifier on the nominal domain ${D}$ such that the following conditions hold:
\begin{assumpttion}[Mollifier]\label{ass:mollifier} The mollifer $\chi(\x):{D} \to [0, 1]$, $\chi \in C^{0, 1} (\overline{{D}})$, has the following properties:
	\begin{enumerate}
		\item[(i)] The mollifier is, on its support, strictly monotonically increasing as a function of $|\x|$ for fixed $\theta(\x)$. 
		\item[(ii)] The mollifier is 0 in a neighborhood of the origin, and 1 on the boundary of ${D}$.
	\end{enumerate}
\end{assumpttion}
Consequently, we define the following transformation satisfying Assumption \ref{ass:phij}:
\begin{definition}[Mollifier transformation]
	With a mollifier satisfying Assumption \ref{ass:mollifier}, we define the mollifer transformation, mapping the nominal domain ${D}$ to the parameterized domain $\mathcal{D}(\y)$, for every $\y\in Y$, by:
	\begin{equation}
		\Phi(\y;\x) = \x + \chi\left(\x\right)\left[ r(\y;\theta(\x)) - r_0(\theta(\x)) \right]\frac{\x}{|\x|},\quad  \x \in D. \label{eq:mollifiertransf}
	\end{equation} 
\end{definition}
From Assumption \ref{ass:mollifier}, we conclude that $\Phi(\y;\x)$  is bijective. Also, it fits the framework \eqref{eq:phidef} with, for $j\geq 1$, 
\begin{equation*}
	\Phi_j(\x):=\chi\left(\x\right)\psi_j(\theta(\x))\frac{\x}{|\x|},\quad \x\in D.
\end{equation*}

To verify the assumptions of Theorem \ref{thm:lpsummability}, we need to calculate an explicit bound on the Jacobian matrices $\D\Phi_j(\x)$, $j \geq 1$, which is given by the following lemma:
\begin{lemma}\label{lem:Dphibound}
For the transformation given by equation \eqref{eq:mollifiertransf} and the parameter-dependent radius as in equation \eqref{eq:rdef}, we have the bound 
\begin{equation}
	|\D\Phi_j(\x)|_{2,2} \leq \chibarone | \psi_j(\theta(\x)) | +  \chibartwo |\psi_j'(\theta(\x))  |, \qquad \text{for a.e. } \x \in {D}, \label{eq:dphijbound}
\end{equation}
with 
\begin{align}
	\chibarone = \sup_{\x \in \mathrm{supp}(\chi)} \frac{1}{\sqrt{2}}\sqrt{ |\nabla_{\x} \chi(\x)|^2  +\frac{\chi(\x)^2}{|\x|^2}  + \sqrt{\left( |\nabla_{\x} \chi|^2  +\frac{\chi(\x)^2}{|\x|^2}   \right)^2-\frac{4\chi(\x)^2}{\r^4}  \left( \nabla_{\x}\chi(\x)\cdot \x \right) ^2}  }
\end{align}
and
\begin{align}
	\chibartwo &= \sup_{\x \in \mathrm{supp}(\chi)}  \frac{\chi\left(\x\right)}{|\x| }.
\end{align}
\end{lemma}
\begin{proof}
	The proof relies on a establishing bound of the form
	\begin{align*}
		|\D\Phi_j(\x)|_{2,2} \leq| \psi_j(\theta(\x)) ||  A_1(\x) |_{2,2}  +  |  \psi_j'(\theta(\x))  || A_2(\x)|_{2,2}, 
	\end{align*}		
	 for a.e. $\x \in D$, and calculating the largest singular values of $A_1(\x)$ and $A_2(\x)$, see Appendix \ref{pf:dphibound} for details.
\end{proof}
Several mollifiers could be used satisfying Assumption \ref{ass:mollifier}, and an important one is the affine mollifier, outlined in the following example:
\begin{example}[Affine mollifier]
	\label{def:linmollifier} 
	\label{col:Dphibound_bound}
	We define the affine mollifier by
	\begin{equation}
		\chi(\x) = \begin{cases}\frac{\x - \frac{r_0^-}{4}}{r_0(\theta(\x))- \frac{r_0^-}{4}}&\text{for }|\x|\geq\frac{r_0^{-}}{4},\\
		0 & \text{otherwise},
\end{cases}		
		\label{eq:mollifier}
	\end{equation}
	where $r_0^-= \inf_{\theta\in [0, 2\pi)} \left(r_0(\theta) \right)$ equals the minimal nominal radius. With the affine mollifier at hand, we can calculate the supremum of $\chi_1(\x)$ and $\chi_2(\x)$ explicitly and obtain, for a.e. $\x \in {D}$,
		\begin{align*}
			\chibarone  = \frac{4}{3{r_0^-}} \sqrt{ \left( 1+ \frac{L_{r_0} ^2}{{r_0^-}^2}    \right) + \left(1 - \frac{r_0^-}{4r_0^-}\right)^2 }, \quad
			\chibartwo = \frac{4}{{3r_0^-}} \left(1 - \frac{r_0^-}{4r_0^-}\right),
		\end{align*}
		where $L_{r_0}$ is the Lipschitz constant of $r_0$.
	For explicit calculations, see Appendix \ref{ap:Dphibound_boundpf}.
\end{example}
Now, we expand on the special case where the functions $\psi_j$ in equation \eqref{eq:rdef} are wavelets. Wavelet expansions allow for localized perturbations in $r(\y; \theta)$, in contrast to the globally supported Fourier modes. For an overview of wavelets, we refer to \cite{Meyer1993}.

To find the optimal convergence rate of the Taylor coefficients \eqref{eq:taylorseq} of the solution to \eqref{eq:transformedvarform}, we propose a sequence $\rho_j$ to verify the assumptions of Theorem \ref{thm:lpsummability}.
Let the wavelets $\{\psi_j(\theta)\}_{j\geq 1}$ be a frame of wavelets generated by periodization of a mother wavelet $\Psi \in C^{0, 1}(\mathbb{R})$ with $\|\Psi\|_{L^\infty(\mathbb{R})}:=1$, 
 as introduced by Yves Meyer \cite{Meyer1993}. As it is more natural, we denote the wavelets by their index $\lambda$, which is a concatenation of their space and scale levels, and using the notation $\abs{\lambda}=l\geq 0$ for the scale level. At each given scale level $l \geq 0$, there are $\mathcal{O}(2^{l})$ wavelets and their supports are such that we can define
\begin{equation}
	M:=\sup_{\theta \in [0,1)} \sum_{k \in \mathbb{Z}}\abs{\Psi(\theta - k)},\quad M_d:=\sup_{\theta \in [0,1)} \sum_{k \in \mathbb{Z}}\abs{\Psi'(\theta - k)}, \label{eq:Mbound}
\end{equation} 
for some $M, M_d \in \mathbb{R}$ independent of the level.
We note that this condition is satisfied by the rapid decay of the mother wavelet \cite{Meyer1993}.
More explicitly, we define, for all $\theta\in[0, 2\pi)$ and $\alpha>2$, 
\begin{align}
	\psi_\lambda(\theta)&:= \frac{r_0^-\vartheta}{M}(1-2^{-\alpha}) 2^{-\alpha \abs{\lambda}} \Psi\left(2^{|\lambda|}(\theta-k)\right), \label{eq:psidef}
\end{align}
where $\vartheta := \sup_{\theta\in [0,2\pi)} \sum_{j \geq 1} |\psi_j(\theta)|$ is the overall maximal amount of shape variation.
Wavelets are such that, due to their localized support, we can sum them very efficiently, as expressed by the following lemma:
\begin{lemma}[Pointwise summability of wavelets] \label{lem:poitwisesum}
	For a Lipschitz-continuous mother wavelet and $\alpha>2$ in \eqref{eq:psidef}, we have the following pointwise bound on the sum of the wavelets and their derivatives:
	\begin{align*}
		\sum_\lambda | \psi_\lambda(\theta) | \leq r_0^-\vartheta, \qquad \sum_\lambda | \psi_\lambda'(\theta ) | \leq r_0^-\vartheta \frac{M_d}{M}\frac{1-2^{-\alpha}}{1-2^{-(\alpha-1)}}, \text{ for a.e. }\theta \in [0,2\pi).
	\end{align*}
	 In particular, the radius expansion fulfills Assumption \ref{ass:phiregularity}.
\end{lemma}
\begin{proof}
	Let $\lambda_l\in \{ \lambda : |\lambda|=l \}$, $l\geq 0$. Using equations \eqref{eq:psidef} and \eqref{eq:Mbound}, we expand the first sum and obtain
	\begin{align}
		\sum_\lambda | \psi_\lambda(\theta) | 
		\leq r_0^-\vartheta(1-2^{-\alpha}) \sum_{l\geq 0} 2^{-\alpha l} 
		\leq r_0^-\vartheta. \label{eq:pointwiseone}
	\end{align}
	Similarly, expanding the second sum yields
		\begin{align}
		\sum_\lambda | \psi_\lambda'(\theta) | 
		\leq r_0^-\vartheta(1-2^{-\alpha}) \frac{M_d}{M} \sum_{l\geq 0} 2^{-(\alpha-1) l} 
		\leq r_0^-\vartheta \frac{M_d}{M}\frac{1-2^{-\alpha}}{1-2^{-(\alpha-1)}}.\label{eq:pointwisetwo}
	\end{align}
\end{proof}
The previous lemma, together with Assumption \ref{ass:mollifier}, ensure that the domain mapping fulfills Assumption \ref{ass:phij}. We turn our attention to the variational formulation on the reference domain, equation \eqref{eq:transformedvarform}.
We remind that Assumption \ref{ass:phij}, when considered with Assumption \ref{ass:a}, ensure the fulfillment of Assumption \ref{ass:afanalytic}.
 the aforementioned properties of the mapping guarantee that, for every $\y\in Y$, $A(\y; \x)$ and $F(\y;\x)$ as defined in \eqref{eq:hatAdef}-\eqref{eq:hatfdef} are in $W^{1,\infty}(D)$ and $L^2(D)$, respectively, with $\y$-independent norm bounds.

To verify the assumptions of Theorem \ref{thm:lpsummability}, we remind that, on the basis of Corollary \ref{lm:Afbound}, the sequence $(b_j(\x))_{j\geq 1}$ appearing in \eqref{eq:AFder} is given by \eqref{eq:bj_phi}, provided we use $\lambda$ as index. Then, for the sequence $(\rho_{\lambda})_{\lambda}$ we propose, for $\beta < \alpha$, 
\begin{equation}
	\rho_\lambda := 1 + \frac{\mathcal{B}_M(1-2^{\beta - \alpha})2^{(\beta - 1)|\lambda|}}{1-2^{-\alpha}},   \label{eq:rhoprop}
\end{equation}
where
\begin{align}
	\mathcal{B}_M = \frac{g_A^{-1}(\amin\sigma_{min}^4)}{1 + \sqrt{2}\chibarone + \sqrt{2}\frac{M_d}{M}\chibartwo} \frac{1}{r_0^- \vartheta}  - 1. \label{eq:B}
\end{align}
We need to ensure that $\mathcal{B}_M>0$ by choosing a small enough value for the maximal shape variation $\vartheta$. In \eqref{eq:B}, we remind the reader that $\sigma_{min}$ is the lower bound on the singular values of the Jacobian matrix. 

To check that \eqref{eq:Kbassumption} is fulfilled, we expand the sum in \eqref{eq:phidef} and treat the resulting terms separately:
\begin{align*}
	\sum_\lambda \rho_\lambda|\D\Phi_\lambda(\x)|_{2,2} 
	&= \sum_\lambda |\D\Phi_\lambda(\x)|_{2,2} + \sum_\lambda \frac{\mathcal{B}_M(1-2^{\beta - \alpha})2^{(\beta - 1)|\lambda|}}{(1-2^{-\alpha})} |\D\Phi_\lambda(\x)|_{2,2}.
\end{align*}
First, making use of Lemmas \ref{lem:Dphibound} and \ref{lem:poitwisesum}, we bound the first term
\begin{align*}
	\sum_\lambda  |\D\Phi_\lambda(\x)|_{2,2}\leq  \sum_\lambda \left(\bar{\chi}_1| \psi_j(\theta(\x)) | + \bar{\chi}_2|  \psi_j'(\theta(\x))  |\right)\leq  r_0^-\vartheta \left( \bar{\chi}_1 + \bar{\chi}_2\frac{M_d}{M} \frac{1-2^{-\alpha}}{1-2^{-(\alpha-1)}}    \right),
\end{align*}
and the second term
\begin{align*}
	&\sum_\lambda   \frac{\mathcal{B}_M(1-2^{\beta - \alpha})2^{(\beta - 1)|\lambda|}}{   (1-2^{-\alpha})} |\D\Phi_\lambda(\x)|_{2,2}\\
	&\qquad\leq\sum_\lambda   \frac{\mathcal{B}_M(1-2^{\beta - \alpha})2^{(\beta - 1)|\lambda|}}{ (1-2^{-\alpha})}\left( \bar{\chi}_1| \psi_\lambda(\theta(\x)) | + \bar{\chi}_2|  \psi_\lambda'(\theta(\x))  |\right)\\
	&\qquad\leq\sum_{l \geq 0} r_0^-\mathcal{B}_M(1-2^{\beta - \alpha})2^{(\beta - 1)l} \left(\bar{\chi}_1 2^{-l} +\bar{\chi}_2 \frac{M_d}{M}\right) 2^{-(\alpha-1) l}
	\leq r_0^- \vartheta  \mathcal{B}_M \left(\bar{\chi}_1 + \bar{\chi}_2\frac{M_d}{M}\right).
\end{align*}
Combining both, we get
\begin{align*}
	\sum_\lambda \rho_\lambda|\D\Phi_\lambda(\x)|_{2,2} &\leq r_0^-\vartheta \left( \bar{\chi}_1 + \bar{\chi}_2\frac{M_d}{M} \frac{1-2^{-\alpha}}{1-2^{-(\alpha-1)}}  + \mathcal{B}_M\left(\bar{\chi}_1 + \bar{\chi}_2\frac{M_d}{M}\right)  \right),
\end{align*}
see also \eqref{eq:pointwiseone} and \eqref{eq:pointwisetwo}.
 In a similar manner, we bound $\sum_\lambda \rho_\lambda|\Phi_\lambda|$ so that, altogether, we obtain 
\begin{align}
	\sum_{\lambda} \rho_\lambda b_\lambda(\x) &\leq r_0^-\vartheta \left( 1+ \sqrt{2} \bar{\chi}_1 +  \sqrt{2}\bar{\chi}_2\frac{M_d}{M} \frac{1-2^{-\alpha}}{1-2^{-(\alpha-1)}} + \mathcal{B}_M \left(1+ \sqrt{2}\chibarone + \sqrt{2}\chibartwo \frac{M_d}{M}\right)  \right) := \Kb, \label{eq:endmollifiercalculation}
\end{align}
for $\vartheta$ small enough to satisfy \eqref{eq:rhoArhoFass}. 
Because of our choice of $\mathcal{B}_M$ as from \eqref{eq:B}, we have
\begin{align*}
	\Kb < g_A^{-1}(\amin\sigma_{min}^4),
\end{align*}
so that equation \eqref{eq:gAdefshape}, and therefore, \eqref{eq:assKbKc} are satisfied. We also have
\begin{align*}
	\rho_\lambda \sim 2^{(\beta-1)|\lambda|},
\end{align*}
which can be reordered \cite{Bachmayr2017a} to obtain 
\begin{align*}
	\rho_j \sim j^{(\beta-1)}, \qquad j  \geq 1.
\end{align*}
This sequence satisfies $(\rho_j^{-1})_{j\geq 1} \, \in  \, \ell^q(\mathbb{N})$ with $q>(\beta-1)^{-1}$, so that we can use Theorem \ref{thm:lpsummability} to obtain $(\anorm{t_\mu})_{\mu \in \mathcal{F}} \, \in \, \ell^p(\mathcal{F})$ for $p>\frac{2(\beta-1)^{-1}}{2+(\beta-1)^{-1}}=\frac{2}{2(\beta-1) + 1}$. This implies that $(\anorm{t_\mu})_{\mu \in \mathcal{F}}$ converges with rate $s=\frac{1}{p} = \frac{2(\beta-1) + 1}{2} = (\beta-1) + \frac{1}{2} = \beta - \frac{1}{2}$, for any $\beta < \alpha$ and $\alpha > 2$.
Therefore, we have the following corollary:
\begin{corollary}\label{cor:mollifiersummability}
	For the solution to the parameterized PDE \eqref{eq:transformedvarform} satisfying Assumption \ref{ass:a}, and a mollifier mapping \eqref{eq:mollifiertransf} with a mollifier satisfying Assumption \ref{ass:mollifier}, and using wavelets as in \eqref{eq:psidef} for the radius expansion \eqref{eq:rdef} with maximal shape variation $\vartheta$ such that $\mathcal{B}_M>0$ and $\Kt < \max\{\rho_A\rho_F,1\}$, we have that
	$(\anorm{t_\mu})_{\mu \in \mathcal{F}}$ $\in \ell^p(\mathcal{F})$, with rate $\frac{1}{p} < \alpha - \frac{1}{2}$, $\alpha > 2$.
\end{corollary}
\begin{remark}[General support basis functions]
	With a technique similar to Remark 2.4 in \cite{Bachmayr2017a}), we can show that, for general support basis functions, we cannot improve over previous results \cite{Chkifa2015,Hiptmair2018}, where a decay rate of $s < \alpha - 1$ was shown.
\end{remark}

\subsection{Harmonic transformation}
\label{sec:harmonic}
Next to the mollifer mapping discussed in the previous subsection, we introduce the harmonic extension to map the star-shaped $D$ to $\mathcal{D}(\y)$, defined implicitly through a Laplace equation on a per spatial dimension basis. For spatial dimension $d=2$, it is defined as the solution to:
\begin{align}
	\begin{cases}
		\Delta \Phi^\kappa = 0 &\text{ in }D,\\
		\Phi^\kappa = \phi^\kappa(\theta)  &\text{ on }\partial D,
	\end{cases} \label{eq:harmonicext}
\end{align}
where $\phi^\kappa(\theta) := f_\kappa(\theta) r(\y; \theta)$ for $\kappa \in \{x, y\}$ and $f_x(\cdot) = \cos(\cdot), f_y(\cdot)=\sin(\cdot)$.
By linearity of the Laplace equation, we can conclude that, for $r(\y;\theta)$ given by \eqref{eq:rdef}, the mapping $\Phi$ satisfies the decomposition \eqref{eq:phidef}.
Moreover, we can conclude that each partial mapping $\Phi_j=[\Phi_j^x, \Phi_j^y]^\top$ solves the harmonic equation with $\Phi_j^\kappa = f_\kappa(\theta) \psi_j(\theta)$ as Dirichlet boundary condition with $\phi_j^\kappa$ defined accordingly.

To show that we can use the harmonic mapping together with the bounds from Corollary \ref{col:Dphibound_bound} and satisfy the assumptions of Theorem \ref{thm:lpsummability}, we will focus our attention to the circular domain $D=\{ \x \in \mathbb{R}^2 : |\x| \leq r_0\}$, with $r_0\in\mathbb{R}_{> 0}$ in the remainder of this section. 
Similar to the mollifier transformation, and in view of \eqref{eq:bj_phi}, we look for pointwise bounds on the Jacobian matrices.

To bound the norm of the Jacobian matrix of each partial transformation $\Phi_j$, we proceed in two step: (I) we bound the supremum of the Jacobian matrix by considering its Value at the boundary only, applying the maximum principle, and (II) we exploit the Dirichlet-to-Neumann map to estimate the derivatives at the boundary in the tangential and normal directions. 

\begin{itemize}
	\item[(I)]
		To bound the norm of the Jacobian matrix of the partial transformation, we bound it, for a.e. $\x \in D$, by the pointwise 1-norm of the elementwise transformations:
		\begin{align*}
			|\D\Phi_j(\x)|_{2,2} 
			\leq\sum_{\kappa \in \{x, y\}}  |\nabla \Phi_j^\kappa(\x)|_1.
		\end{align*}
		As $\Phi_j$ is a weakly harmonic function, the derivatives of $\Phi_j$ are weakly harmonic as well, and the absolute value of the derivatives of $\Phi_j$ and their sum are weakly subharmonic. 
		Hence, we can apply the maximum principle \cite{Segala1999} to conclude that:
		\begin{align}
			\sup_{\x \in D} \sum_{\kappa \in \{x, y\}}  |\nabla \Phi_j^\kappa (\x)|_1 \leq \sup_{\x \in \partial D} \sum_{\kappa \in \{x, y\}}  |\nabla \Phi_j^\kappa(\x)|_1. \label{eq:partialmaxprinc}
		\end{align}				
			\item[(II)]
		To evaluate the right-hand side of \eqref{eq:partialmaxprinc} at the boundary, we estimate, for a.e. $\x \in D$,
		\begin{align*}
			|\nabla \Phi_j^\kappa(\x)|_1 \leq  \sqrt{2} |\nabla \Phi_j^\kappa(\x)|,
		\end{align*}
		where we remind that $|\cdot|$ is the Euclidean norm. Now, we calculate the derivatives explicitly at the boundary in both the radial and the angular directions, that is, respectively the normal and tangential derivatives at the boundary. 
		For the angular derivative, we can just differentiate the boundary condition to obtain $ \frac{\partial \phi_j^\kappa}{\partial \theta} = (f_\kappa\psi_j)'(\theta)$, for $\theta \in [0,2\pi)$ and $\kappa\in\{x,y\}$. 
		To complete the computation of the gradient, we observe that the radial derivative of $\Phi_j^\kappa$ equals to the Dirichlet-to-Neumann (DtN) map of $\phi_j^\kappa$. 
		For circular domains, the DtN map of the Laplacian is given by $\sqrt{-\Delta_{\partial D}}$, the square root of the Laplace-Beltrami operator on the boundary circle  \cite{Girouard2022}.
		Because our domain is circular, we pick up a factor of $\frac{1}{r_0}$ due to the curvature of the circle, and we can expand using the symbol of the operator and the Fourier transform $\mathcal{F}(\cdot)$,
		\begin{align*}
			\sqrt{-\Delta_{\partial D}} \phi_j^\kappa  
			&=\frac{1}{r_0}\int_\mathbb{R} |\xi| \mathcal{F}(\phi_j^\kappa )(\xi) e^{-i \xi x} \d\xi
			=\frac{1}{r_0}\int_\mathbb{R} i\xi (-i) \sign(\xi) \mathcal{F}(\phi_j^\kappa )(\xi) e^{-i \xi x} \d\xi\\
			&=\frac{1}{r_0}\int_\mathbb{R} \mathcal{F}\left(\frac{\d}{\d x}\mathcal{H}(\phi_j^\kappa )\right) (\xi) e^{-i \xi x} \d\xi
			=\frac{1}{r_0}\frac{\d}{\d x}\mathcal{H}(\phi_j^\kappa ),
		\end{align*}
		where $\mathcal{H}(\cdot)$ denotes the Hilbert transform \cite{Hahn1996}.

		Now, we remember that for any signal $f$, the signal $f^A:=f+i\mathcal{H}(f)$ is analytic, where analyticity is meant in the signal processing sense; it has no negative frequency components. Moreover, we recall that the instantaneous amplitude $A(f^A)(x)$ of any analytic signal $f^A$ is given by the modulus of the signal itself. Therefore we observe that, for every $\theta\in[0,2\pi)$ and $\kappa \in \{x,y\}$,
		\begin{align}
			\left| \nabla \phi_j^\kappa(\theta)  \right|&\leq \sqrt{2} \max\{1, r_0^{-1}\} \left|(f_\kappa\psi_j)'(\theta) + i \mathcal{H}\left( (f_\kappa\psi_j)' \right)(\theta)  \right| \nonumber\\
			&= \sqrt{2} \max\{1, r_0^{-1}\} A((f_\kappa\psi_j)')(\theta) \nonumber\\
			&\leq \sqrt{2} \max\{1, r_0^{-1}\} \left( A(\psi_j)(\theta) + A(\psi_j')(\theta) \right),
		\end{align}
		where we have used the definition of $f_\kappa$ which have amplitude 1, and their derivatives as well. From this, for a.e. $\theta \in [0,2\pi)$, we obtain the bound 
	\begin{align*}
		|\D\Phi_j(\theta)|_{2,2} \leq  4 \max\{1, r_0^{-1}\}  \left( A(\psi_j)(\theta) + A(\psi_j')(\theta) \right).
		\end{align*}

\end{itemize}
Moreover, the sum $\sum_{j \geq 1} \rho_j|\D\Phi_j(\x)|_{2,2} $ is bounded from above by the weakly subharmonic function\\ $\sum_{j\geq 1}  \sum_{\kappa \in \{x, y\}}  |\nabla \phi_j^\kappa(\x)|_1$ and its maximum is attained at the boundary as well. Therefore, we can repeat steps (I) and (II) to obtain the bound
\begin{align}
	\sup_{\x \in D}\sum_{j\geq 1} \rho_j |\D\Phi_j(\x)|_{2,2}  &\leq \sup_{\theta \in [0,2\pi)}  4 \max\{1, r_0^{-1}\} \left(  \sum_{j \geq 1}\rho_j A(\psi_j)(\theta) + \sum_{j \geq 1}\rho_j A(\psi_j')(\theta) \right). \label{eq:rhodphisum}
\end{align}
Similarly, since $\Phi_j$ is harmonic, we can repeat the same argument and obtain 
\begin{align}
	\sup_{\x \in D}\sum_{j\geq 1} \rho_j |\Phi_j(\x)|\leq \sup_{\theta \in [0,2\pi)} \sum_j \rho_j |\psi_j(\theta)|.\label{eq:rhophisum}
\end{align}
Finally, combining equations \eqref{eq:rhodphisum} and \eqref{eq:rhophisum}, we arrive at the following lemma:
\begin{lemma}\label{lem:harmboundbj}
	For the harmonic transformation given by equation \eqref{eq:harmonicext}, with $b_j(\x)$ from equation \eqref{eq:bj_phi},  and a positive sequence $(\rho_j)_{j \geq 1}$, on a circular, two-dimensional reference domain, we have the bound 
	\begin{align*}
		\sup_{\x \in D}&\sum_j \rho_j b_j(\x) \leq \\ &\sup_{\theta \in [0,2\pi)} \left( \left(4\sqrt{2} \max\{1, r_0^{-1}\} + 1\right) \sum_{j}\rho_j A(\psi_j)(\theta) +  4\sqrt{2} \max\{1, r_0^{-1}\} \sum_{j}\rho_j A(\psi_j')(\theta) \right).
	\end{align*}
\end{lemma}
Similarly to Subsection \ref{sec:explicittransformation}, we will focus on the wavelet case, with $\alpha>2$, Lipschitz mother wavelet, and we index the wavelets using $\lambda$ instead of $j$. To show the desired estimate on $\sum_{\lambda} \rho_\lambda b_\lambda (\x)$, we set, similarly to equation \eqref{eq:Mbound}:
\begin{equation}
	M_H:=\sup_{\theta \in [0,1)} \sum_{k \in \mathbb{Z}}A(\Psi)(\theta - k),\quad M_{d,H}:=\sup_{\theta \in [0,1)} \sum_{k \in \mathbb{Z}}A(\Psi')(\theta - k), \label{eq:Mbound2}
\end{equation} 
and we define $\psi_\lambda$ similarly to \eqref{eq:psidef}, and proceed as in Subsection \ref{sec:explicittransformation}, replacing $M$ and $M_d$ with $M_H$, $M_{d,H}$ whenever applicable. We have to be careful though, as the Hilbert transform does not preserve the locality properties of $\Psi$. However, if the mother wavelet $\Psi$ has $n$ vanishing moments, we have $\left|  \Psi' (\theta) \right| \leq \mathcal{O}(n^{-1-n})$ \cite{Chaudhury2011}, which gives us the required summability when we periodize to obtain a wavelet system on the circle.  Hence, we recover Lemma \ref{lem:poitwisesum} with the aforementioned replacements.

Now, we take $\rho_\lambda$  as defined in equation \eqref{eq:rhoprop}, where $\mathcal{B}_M$ is replaced by 
\begin{equation}
	\mathcal{B}_H = \frac { \frac{1}{\vartheta}g_A^{-1}(\amin\sigma_{min}^4) - 4\sqrt{2}\max\{1, r_0 \}   - \left(4\sqrt{2}\max\{1, r_0 \} + r_0 \right)   \frac{M_{d,H}}{M_H}\frac{1-2^{-\alpha}}{1-2^{-(\alpha - 1)}}  }{4\sqrt{2}\max\{1, r_0 \}\left( 1+\frac{M_{d,H}}{M_H} \right) + r_0  }.
\end{equation}
Similarly to $\mathcal{B}_M$, we should be careful in choosing $\vartheta$ small enough such that we keep $\mathcal{B}_H$ positive.
By repeating the calculations from equation \eqref{eq:rhoprop} to Corollary \ref{cor:mollifiersummability}, we arrive at the following statement:
\begin{corollary}\label{cor:harmonicsummability}
	For the solution to the parameterized PDE \eqref{eq:transformedvarform} satisfying Assumption \ref{ass:a}, and a harmonic mapping \eqref{eq:harmonicext}, and using wavelets as in \eqref{eq:psidef} for the radius expansion \eqref{eq:rdef} with maximal shape variation $\vartheta$ such that $\mathcal{B}_H>0$ and $\Kt < \max\{\rho_A\rho_F,1\}$, we have that
	$(\anorm{t_\mu})_{\mu \in \mathcal{F}}$ $\in \ell^p(\mathcal{F})$, with rate $\frac{1}{p} < \alpha - \frac{1}{2}$,  $\alpha > 2$.
\end{corollary}

\begin{remark} \label{rem:sepofvar}
	The results leading up to Lemma \ref{lem:harmboundbj} that have been obtained using the symbol of $\sqrt{-\Delta_{\partial D}}$ can be obtained via direct computations of the exact solution of equation \eqref{eq:harmonicext} which, in turn, can be obtained via separation of variables with Lipschitz boundary, see for instance \cite{Natalini2008}.
\end{remark}

\begin{remark}
	For non circular star-shaped domains with Lipschitz boundary, we cannot rely on an explicit formula for the DtN map similarly to the approach outlined in this section, see for instance \cite{Girouard2022}. 
	Direct computations as outlined in Remark \ref{rem:sepofvar} do not seem to help drawing useful conclusions either. We do believe that the general result should follow from quasi-local properties of the DtN map, but these are hard to find in the literature which focuses more on global properties expressed in terms of norms.
\end{remark}

%% file: 5numerical_results/5numerical_results.tex
\section{Numerical illustrations}
\label{sec:numerical}
To illustrate the theoretical results obtained in Section \ref{sec:modelproblem}, we compute the Taylor coefficients of the solution to the parameterized shape problem \eqref{eq:transformedvarform} and estimate their decay. In equation \eqref{eq:transformedvarform}, we use $a(\bm{x}) = f(\bm{x}) \equiv 1$, and we take the circle with $r_0\equiv 1$ as the star-shaped reference domain.

To calculate the Taylor coefficients, we have implemented the Alternating greedy Taylor algorithm introduced in \cite{Cohen2015}. The algorithm iteratively builds up a downward closed index set $\Lambda$ by determining its reduced set of neighbors $\mathcal{N}(\Lambda) \subset \mathcal{F}$ \cite{Cohen2015}: 
\begin{equation*}
	\mathcal{N}(\Lambda):= \left\{   \nu \notin \Lambda \text{ such that } \Lambda \cup \{\nu\} \text{ is downward closed and } \supp(\nu) \leq \max_{\mu \in \Lambda}\left\{\supp(\mu)\right\} +1 \right\},
\end{equation*}
and adding alternatingly the largest Taylor coefficient and the Taylor coefficient that has been in the reduced set of neighbors for the largest number of iterations. This way, all Taylor coefficients would be added to $\Lambda$ \cite{Cohen2015} if we did not put a limit on the number of iterations. In each test, we continue until we have added 1000 coefficients to $\Lambda$. To estimate the decay rate, we follow \cite{Bachmayr2017a} and introduce the decreasing rearrangement $(t^*_n)_{n\geq 1}$ of $(\|t_\mu\|_{V})_{\mu \in \mathcal{F}}$. Then, $(\|t_\mu\|_{V})_{\mu \in \mathcal{F}} \in \ell^p(\mathcal{F})$ implies that, for some $C>0$, we have $t^*_n \leq Cn^{-\frac{1}{p}}$, and, if $t^*_n \leq Cn^{-\frac{1}{q}}$ for some $C,q > 0$, we have that $(\|t_\mu\|_{V})_{\mu \in \mathcal{F}} \in \ell^p(\mathcal{F})$ for any $p > q$. Hence, we can estimate the limiting decay rate of $(\|t_\mu\|_{V})_{\mu \in \mathcal{F}}$ from the largest $s$ such that $\sup_{n \geq 1} n^st^*_n$ is finite. To estimate this largest $s$, we perform a linear least squares regression on an exponentially distributed sample of $n$  to find the decay rate in a log-log scale. By sampling exponentially from the decreasing rearrangement $(t^*_n)_{n \geq 1}$, we mitigate the staircasing effect introduced by the wavelet expansion. To make sure we do not include any pre-asymptotic phase, we skip the largest coefficients in the pre-asymptotic behavior and, due to the alternating nature of the alternating greedy algorithm, we exclude the smallest half of the coefficients, because these coefficients may be included due to the alternating steps including the oldest coefficients instead of the largest ones. 

We have implemented the alternating greedy algorithm in Python. To compute the Taylor coefficients, we are using the finite element library Dolfinx \cite{Alnes2014,Scroggs2022a,Scroggs2022}, employing the linear system solver PETSc \cite{Brown2022,Dalcin2011} and the domain meshing library Gmsh \cite{Geuzaine2009} to generate the triangular, unstructured reference mesh. To evaluate the wavelets and scaling functions efficiently, we have used the PyWavelet library \cite{Lee2019}. Since our theoretical results do not take the finite element error in the computation of the Taylor coefficients into account, we have, apart from Subsection \ref{sec:convstudy}, chosen the mesh to be fine enough so that its effect on the decay of the Taylor coefficients is negligible. To take the radial dependence of the mollifier mapping \eqref{eq:mollifiertransf} into account, we are computing the Taylor coefficients using a mesh with a mesh size decreasing linearly with the distance from the origin. When using the harmonic mapping instead, we have used a mesh with a mesh size proportional to the inverse of the distance to the boundary. Both meshes, when using either the mollifier or harmonic mapping, have been scaled so as to have identical spatial discretization at their boundaries. 

However, we mention that, although we have shown a possible convergence rate in the previous sections and we observe good approximation of the expected decay in practice, we cannot prove theoretically to attain this maximal convergence rate using the alternating greedy algorithm \cite{Cohen2015}.

When calculating the required Taylor coefficients, we need to evaluate the derivative $\partial^\mu A(\bm{0})$ at each step. Due to the multiplicative structure of $A(\y)$, see \eqref{eq:hatAdef}, the computational size of $\partial^\mu A(\bm{0})$  scales super-exponentially with the order $\mu$, and, for memory reasons, we had to truncate $\mathcal{F}$ and discard any multi-index $\mu$ with $|\mu| > 7$. For all computations with maximal shape variation $\vartheta < 32\%$, this does not impact the downwards closed index set $\Lambda$, as the corresponding Taylor coefficients do not enter $\Lambda$. For $\vartheta \geq 32\%$, the number of skipped coefficients is very small, and therefore, this does not impact the obtained convergence rates either. 

To verify the theoretical convergence rates when using \eqref{eq:rdef}, we define the Fourier basis $\{\psi^F_j\}_{j\geq 1}$ by
\begin{align*}
	\psi^F_j(\theta) \sim \begin{cases}
		\left(\frac{j}{2}\right)^{-\alpha} \sin(\frac{j}{2}\theta) &\text{for } j\text{ even},\\
		\left(\frac{j+1}{2}\right)^{-\alpha} \cos(\frac{j+1}{2}\theta) & \text{for } j\text{ odd},
	\end{cases}
\end{align*}
where $\sim$ entails the constant we use to rescale the Fourier basis such that the total shape variation  equals $\vartheta$. Next to the Fourier expansion, we use a wavelet basis by periodization \cite{Meyer1993} of the orthogonal Daubechies4 \cite{Daubechies1992} wavelet, obtained via equation \eqref{eq:psidef}.  In both cases, the series expansion in \eqref{eq:rdef} is scaled to have maximal shape variation  $\vartheta$.


In the remainder of this section, we present the numerical verification of our analytical results obtained in Subsection \ref{sec:application}. In Subsection \ref{sec:vartheta5percent}, we present the numerical rates for fixed maximal shape variations of $\vartheta = 5\%$, in Subsection \ref{sec:shapevar} we discuss the effect of the maximal shape variations $\vartheta$ on the rates, and finally, in Subsection \ref{sec:convstudy}, we show the relation between the spatial discretization and the observed rates.

\begin{figure}[t!]
\ifgraphs
\centering
	\begin{subfigure}[t]{0.3\textwidth}
		\begin{tikzpicture}[scale=0.75]
	   	    \begin{axis}[
	   	    legend style={at={(0.02,0.02)},anchor=south west, nodes={scale=0.8, transform shape}},
	   	    height =12cm,
	   	    width = 1.5\linewidth, 
	   	    xlabel = {$n$},
	        ylabel = {\ylab},
	        grid=both,
	        major grid style={black!50},
	        xmode=log, ymode=log,
	        xmin=1e0, xmax=5e2,
	        ymin=1e-11, ymax=1e-0,
	        yticklabel style={
	            /pgf/number format/fixed,
	            /pgf/number format/precision=0
	        },
	        scaled y ticks=false,
	        xticklabel style={
	            /pgf/number format/fixed,
	            /pgf/number format/precision=0
	        },
	        scaled y ticks=false,
	        legend columns=2
	        title={$\alpha = 3$},
	        title style={
	        	font=\LARGE
	        },
	        ]

	   	    \plotfit{Alternating_greedy_fourier_alpha_2_h1_theta005_mollifier_deg2_sourceconst_final}{color_llrr}{M--F\hspace*{0.1cm}}
	   	    \plotfit{Alternating_greedy_daubechies4_alpha_2_h1_theta005_mollifier_deg2_sourceconst_final}{color_llrrg}{M--W\hspace*{0.1cm}}
	   	    \plotfit{Alternating_greedy_fourier_alpha_2_h1_theta005_harmonic_deg2_sourceconst_final}{color_llgg}{H--F\hspace*{0.1cm}}
	   	    \plotfit{Alternating_greedy_daubechies4_alpha_2_h1_theta005_harmonic_deg2_sourceconst_final}{color_llbb}{H--W\hspace*{0.1cm}}
	   	    
        	\end{axis}
		\end{tikzpicture}
		\caption{$\alpha=2$}
	\end{subfigure}
	~
	\begin{subfigure}[t]{0.3\textwidth}
		\begin{tikzpicture}[scale=0.75]
	   	    \begin{axis}[
	   	    legend style={at={(0.02,0.02)},anchor=south west, nodes={scale=0.8, transform shape}},
	   	    height =12cm,
	   	    width = 1.5\linewidth, 
	   	    xlabel = {$n$},
	        grid=both,
	        major grid style={black!50},
	        xmode=log, ymode=log,
	        xmin=1e0, xmax=5e2,
	        ymin=1e-11, ymax=1e-0,
	        yticklabel style={
	            /pgf/number format/fixed,
	            /pgf/number format/precision=0,
	            text opacity=0
	        },
	        scaled y ticks=false,
	        xticklabel style={
	            /pgf/number format/fixed,
	            /pgf/number format/precision=0
	        },
	        scaled y ticks=false,
	        legend columns=2
	        ]
	        
	   	    \plotfit{Alternating_greedy_fourier_alpha_25_h1_theta005_mollifier_deg2_sourceconst_final}{color_rr}{M--F\hspace*{0.1cm}}
	   	    \plotfit{Alternating_greedy_daubechies4_alpha_25_h1_theta005_mollifier_deg2_sourceconst_final}{color_rrg}{M--W\hspace*{0.1cm}}
	   	    \plotfit{Alternating_greedy_fourier_alpha_25_h1_theta005_harmonic_deg2_sourceconst_final}{color_gg}{H--F\hspace*{0.1cm}}
	   	    \plotfit{Alternating_greedy_daubechies4_alpha_25_h1_theta005_harmonic_deg2_sourceconst_final}{color_bb}{H--W\hspace*{0.1cm}}
	   	    
	        \end{axis}
		\end{tikzpicture}
		\caption{$\alpha=2.5$}
	\end{subfigure}
	~
	\begin{subfigure}[t]{0.3\textwidth}
		\begin{tikzpicture}[scale=0.75]
	   	    \begin{axis}[
	   	    legend style={at={(0.02,0.02)},anchor=south west, nodes={scale=0.8, transform shape}},
	   	    height =12cm,
	   	    width = 1.5\linewidth, 
	   	    xlabel = {$n$},
	        grid=both,
	        major grid style={black!50},
	        xmode=log, ymode=log,
	        xmin=1e0, xmax=5e2,
	        ymin=1e-11, ymax=1e-0,
	        yticklabel style={
	            /pgf/number format/fixed,
	            /pgf/number format/precision=0,
	            text opacity=0
	        },
	        scaled y ticks=false,
	        xticklabel style={
	            /pgf/number format/fixed,
	            /pgf/number format/precision=0
	        },
	        scaled y ticks=false,
	        legend columns=2
	        ]
	        
	   	    \plotfit{Alternating_greedy_fourier_alpha_3_h1_theta005_mollifier_deg2_sourceconst_final}{color_ddrr}{M--F\hspace*{0.1cm}}
	   	    \plotfit{Alternating_greedy_daubechies4_alpha_3_h1_theta005_mollifier_deg2_sourceconst_final}{color_ddrrg}{M--W\hspace*{0.1cm}}
	   	    \plotfit{Alternating_greedy_fourier_alpha_3_h1_theta005_harmonic_deg2_sourceconst_final}{color_ddgg}{H--F\hspace*{0.1cm}}
	   	    \plotfit{Alternating_greedy_daubechies4_alpha_3_h1_theta005_harmonic_deg2_sourceconst_final}{color_ddbb}{H--W\hspace*{0.1cm}}
	   	    
	        \end{axis}
		\end{tikzpicture}
		\caption{$\alpha=3$}
	\end{subfigure}
\else
	\begin{subfigure}[t]{0.37\textwidth}
		\includegraphics[scale=0.075]{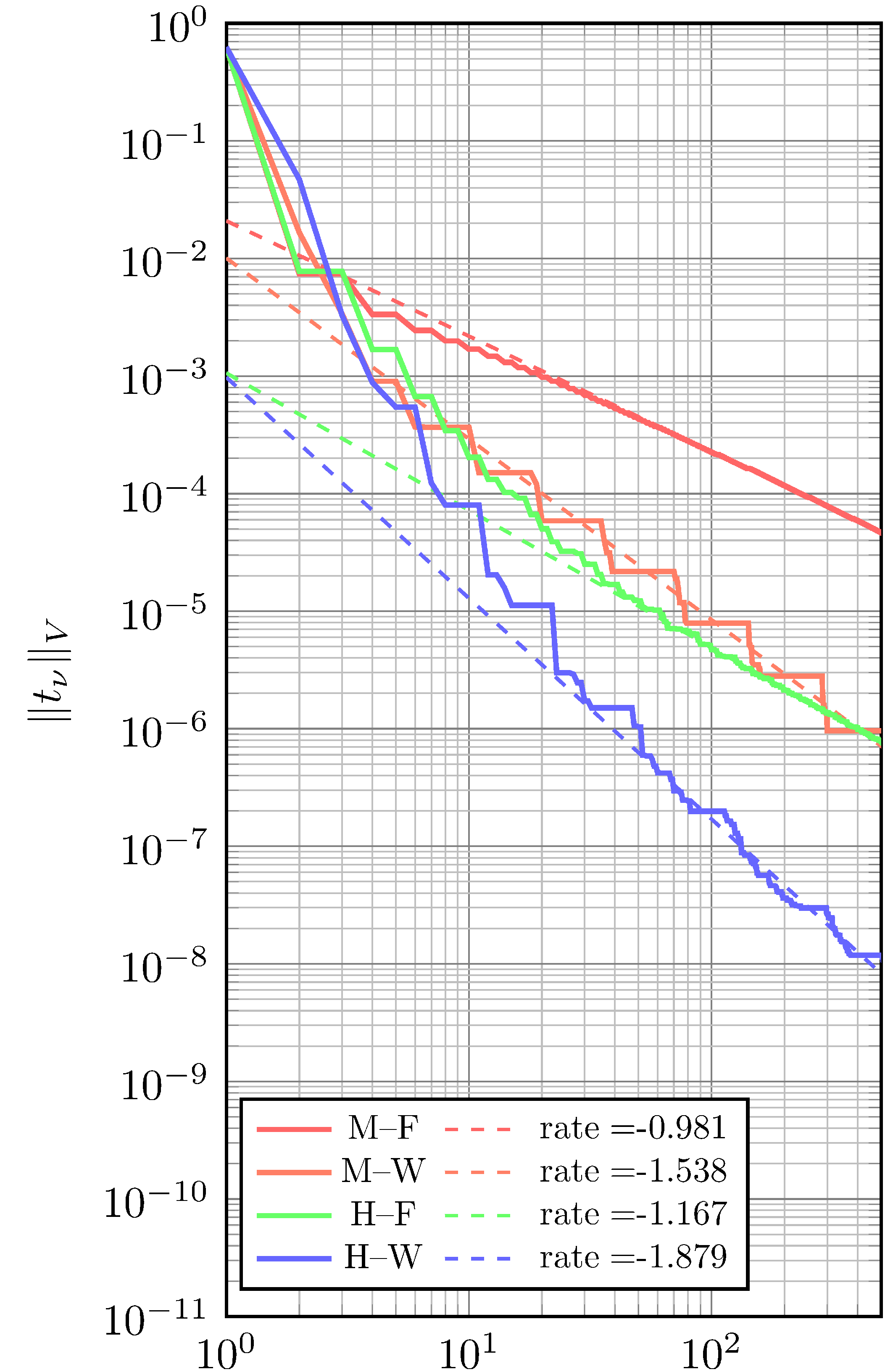}
		\caption{$\alpha=2$}
	\end{subfigure}
	~
	\begin{subfigure}[t]{0.29\textwidth}
		\includegraphics[scale=0.075]{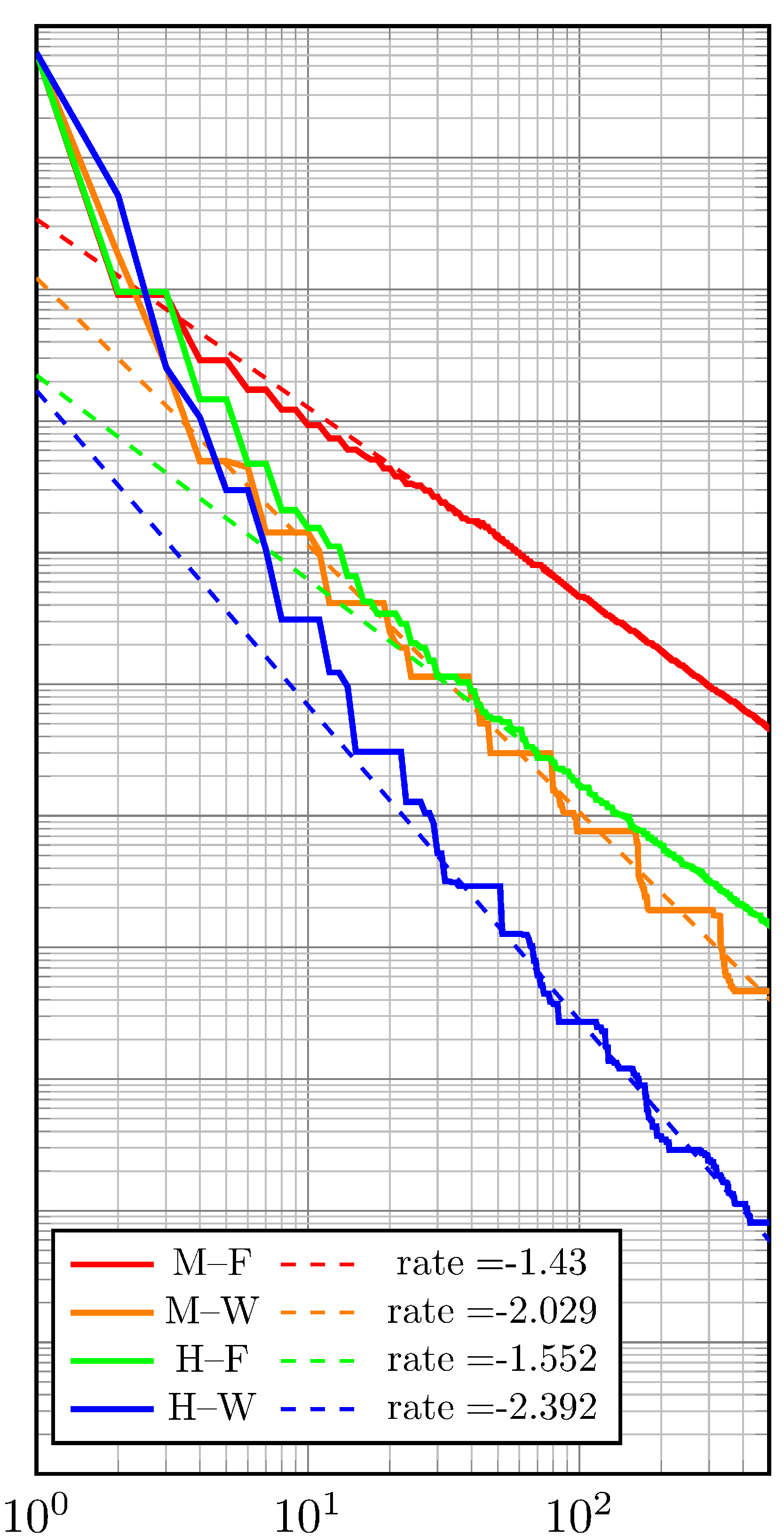}
		\caption{$\alpha=2.5$}
	\end{subfigure}
	~
	\begin{subfigure}[t]{0.29\textwidth}
		\includegraphics[scale=0.075]{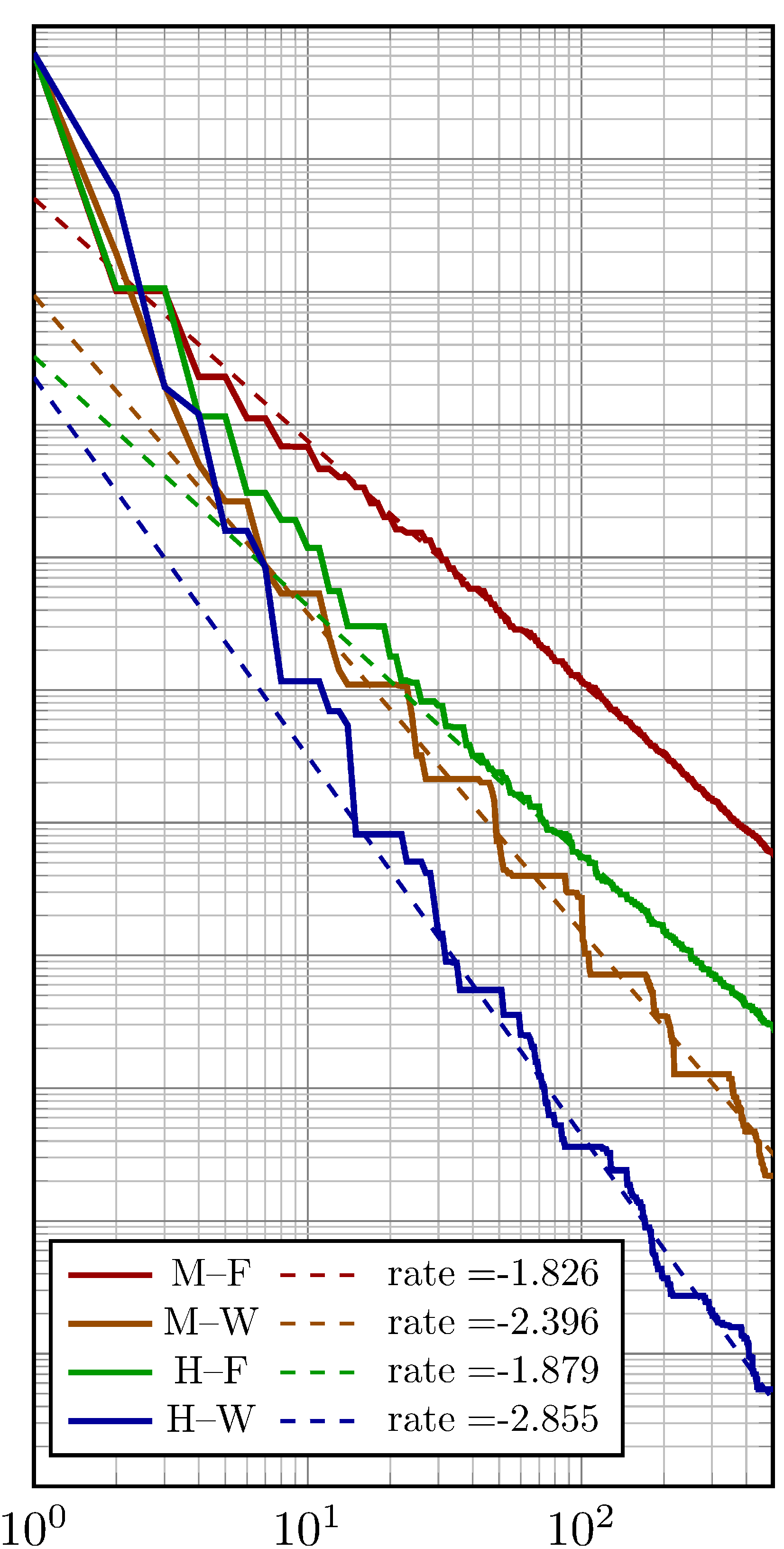}
		\caption{$\alpha=3$}
	\end{subfigure}
\fi
	\caption{Decay of the $V$-norm of the Taylor coefficients of the PDE solution for the Fourier expansion with mollifier mapping (M--F), Wavelet expansion with mollifier mapping (M--W), Fourier expansion with harmonic mapping (H--F), and Wavelet expansion with harmonic mapping (H--W), for $\alpha = 2$ (a), $\alpha=2.5$ (b), and $\alpha=3$ (c), and with $\vartheta=5\%$ maximal shape variations. We denote by $n$ the index in the decreasing rearrangement $(t^*_n)_{n \geq 1}$.}
	\label{fig:threeresults}
\end{figure}

\begin{figure}
	\centering 
	\ifgraphs
	\begin{subfigure}[t]{0.425\textwidth}
		\begin{tikzpicture}[scale=0.75]
	   	    \begin{axis}[
	   	    legend style={at={(0.02,0.02)},anchor=south west, nodes={scale=0.8, transform shape}},
	   	    height =12cm,
	   	    width = 1.45\linewidth, 
	   	    xlabel = {$n$},
	        ylabel = {\ylab},
	        grid=both,
	        major grid style={black!50},
	        xmode=log, ymode=log,
	        xmin=1e0, xmax=5e2,
	        ymin=1e-11, ymax=1e-0,
	        yticklabel style={
	            /pgf/number format/fixed,
	            /pgf/number format/precision=0,
	        },
	        scaled y ticks=false,
	        xticklabel style={
	            /pgf/number format/fixed,
	            /pgf/number format/precision=0
	        },
	        scaled y ticks=false,
	        legend columns=2
	        ]
	        
	   	    \plotfit{Alternating_greedy_fourier_alpha_2_h1_theta005_mollifier_deg2_sourceconst_final}{color_llrr}{$\alpha=2$\hspace*{0.1cm}}
	   	    \plotfit{Alternating_greedy_fourier_alpha_25_h1_theta005_mollifier_deg2_sourceconst_final}{color_rr}{$\alpha=2.5$\hspace*{0.1cm}}
	   	    \plotfit{Alternating_greedy_fourier_alpha_3_h1_theta005_mollifier_deg2_sourceconst_final}{color_ddrr}{$\alpha=3$\hspace*{0.1cm}}
	   	    
	        \end{axis}
		\end{tikzpicture}
		\caption{Mollifier mapping with Fourier expansion}
	\end{subfigure}
	~
	\begin{subfigure}[t]{0.47\textwidth}
		\begin{tikzpicture}[scale=0.75]
	   	    \begin{axis}[
	   	    legend style={at={(0.02,0.02)},anchor=south west, nodes={scale=0.8, transform shape}},
	   	    height =12cm,
	   	    width = 1.3\linewidth, 
	   	    xlabel = {$n$},
	        grid=both,
	        major grid style={black!50},
	        xmode=log, ymode=log,
	        xmin=1e0, xmax=5e2,
	        ymin=1e-11, ymax=1e-0,
	        yticklabel style={
	            /pgf/number format/fixed,
	            /pgf/number format/precision=0,
	            text opacity=0
	        },
	        scaled y ticks=false,
	        xticklabel style={
	            /pgf/number format/fixed,
	            /pgf/number format/precision=0
	        },
	        scaled y ticks=false,
	        legend columns=2
	        ]
	        
	   	    \plotfit{Alternating_greedy_daubechies4_alpha_2_h1_theta005_mollifier_deg2_sourceconst_final}{color_llrrg}{$\alpha=2$\hspace*{0.1cm}}
	   	    \plotfit{Alternating_greedy_daubechies4_alpha_25_h1_theta005_mollifier_deg2_sourceconst_final}{color_rrg}{$\alpha=2.5$\hspace*{0.1cm}}
	   	    \plotfit{Alternating_greedy_daubechies4_alpha_3_h1_theta005_mollifier_deg2_sourceconst_final}{color_ddrrg}{$\alpha=3$\hspace*{0.1cm}}
	   	    
	        \end{axis}
		\end{tikzpicture}
		\caption{Mollifier mapping with wavelet expansion}
	\end{subfigure}\\
	\begin{subfigure}[t]{0.425\textwidth}
		\begin{tikzpicture}[scale=0.75]
	   	    \begin{axis}[
	   	    legend style={at={(0.02,0.02)},anchor=south west, nodes={scale=0.8, transform shape}},
	   	    height =12cm,
	   	    width = 1.45\linewidth, 
	   	    xlabel = {$n$},
	        ylabel = {\ylab},
	        grid=both,
	        major grid style={black!50},
	        xmode=log, ymode=log,
	        xmin=1e0, xmax=5e2,
	        ymin=1e-11, ymax=1e-0,
	        yticklabel style={
	            /pgf/number format/fixed,
	            /pgf/number format/precision=0,
	        },
	        scaled y ticks=false,
	        xticklabel style={
	            /pgf/number format/fixed,
	            /pgf/number format/precision=0
	        },
	        scaled y ticks=false,
	        legend columns=2
	        ]
	        
	   	    \plotfit{Alternating_greedy_fourier_alpha_2_h1_theta005_harmonic_deg2_sourceconst_final}{color_llgg}{$\alpha=2$\hspace*{0.1cm}}
	   	    \plotfit{Alternating_greedy_fourier_alpha_25_h1_theta005_harmonic_deg2_sourceconst_final}{color_gg}{$\alpha=2.5$\hspace*{0.1cm}}
	   	    \plotfit{Alternating_greedy_fourier_alpha_3_h1_theta005_harmonic_deg2_sourceconst_final}{color_ddgg}{$\alpha=3$\hspace*{0.1cm}}
	   	    
	        \end{axis}
		\end{tikzpicture}
		\caption{Harmonic mapping with Fourier expansion}
	\end{subfigure}
	~
	\begin{subfigure}[t]{0.47\textwidth}
		\begin{tikzpicture}[scale=0.75]
	   	    \begin{axis}[
	   	    legend style={at={(0.02,0.02)},anchor=south west, nodes={scale=0.8, transform shape}},
	   	    height =12cm,
	   	    width = 1.3\linewidth, 
	   	    xlabel = {$n$},
	        grid=both,
	        major grid style={black!50},
	        xmode=log, ymode=log,
	        xmin=1e0, xmax=5e2,
	        ymin=1e-11, ymax=1e-0,
	        yticklabel style={
	            /pgf/number format/fixed,
	            /pgf/number format/precision=0,
	            text opacity=0
	        },
	        scaled y ticks=false,
	        xticklabel style={
	            /pgf/number format/fixed,
	            /pgf/number format/precision=0
	        },
	        scaled y ticks=false,
	        legend columns=2
	        ]
	        
	   	    \plotfit{Alternating_greedy_daubechies4_alpha_2_h1_theta005_harmonic_deg2_sourceconst_final}{color_llbb}{$\alpha=2$\hspace*{0.1cm}}
	   	    \plotfit{Alternating_greedy_daubechies4_alpha_25_h1_theta005_harmonic_deg2_sourceconst_final}{color_bb}{$\alpha=2.5$\hspace*{0.1cm}}
	   	    \plotfit{Alternating_greedy_daubechies4_alpha_3_h1_theta005_harmonic_deg2_sourceconst_final}{color_ddbb}{$\alpha=3$\hspace*{0.1cm}}
	   	    
	        \end{axis}
		\end{tikzpicture}
		\caption{Harmonic mapping with wavelet expansion}
	\end{subfigure}
	\else
	\begin{subfigure}[t]{0.49\textwidth}
		\includegraphics[scale=0.074]{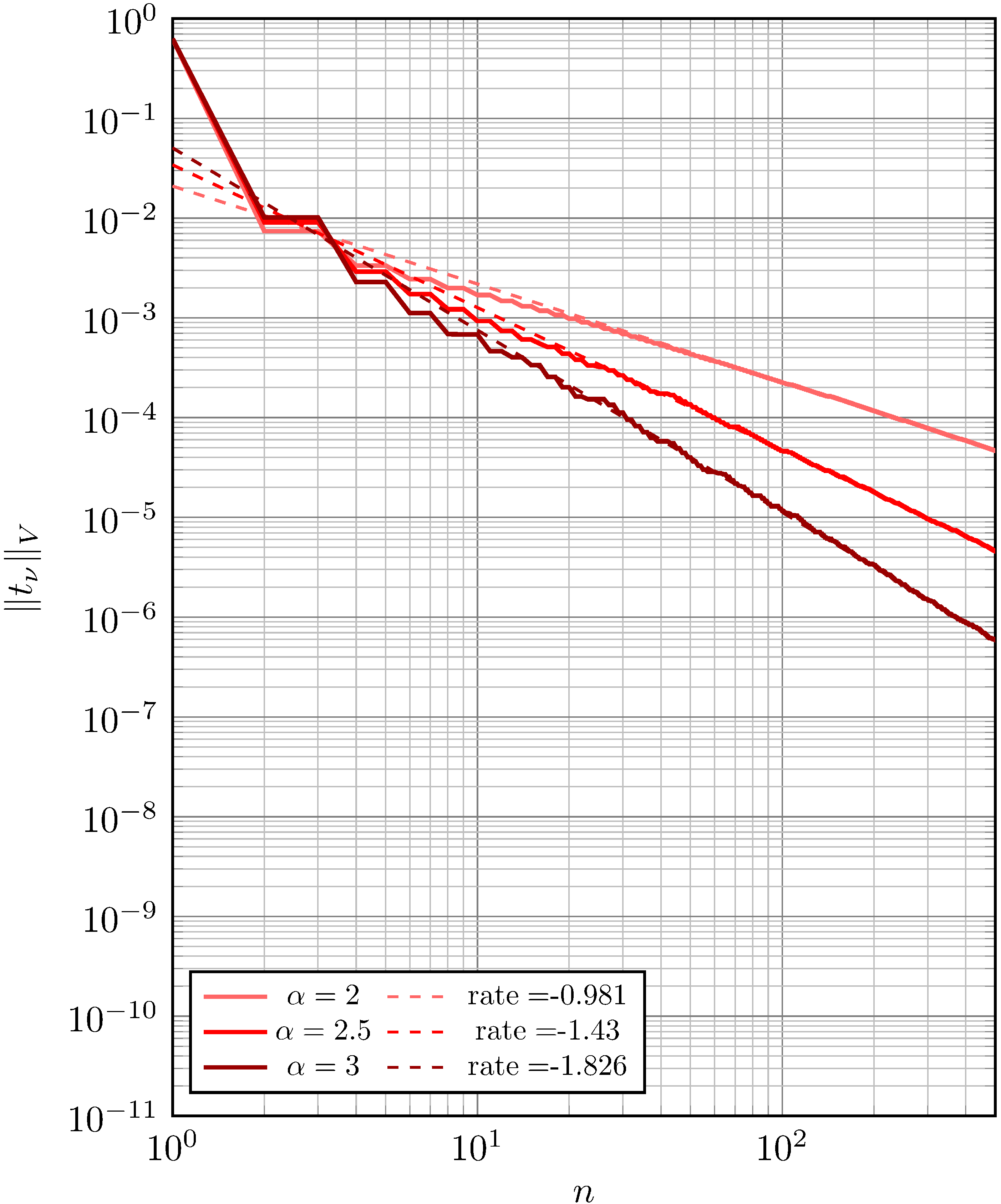}
		\caption{Mollifier mapping with Fourier expansion}
	\end{subfigure}
	~
	\begin{subfigure}[t]{0.45\textwidth}
		\includegraphics[scale=0.074]{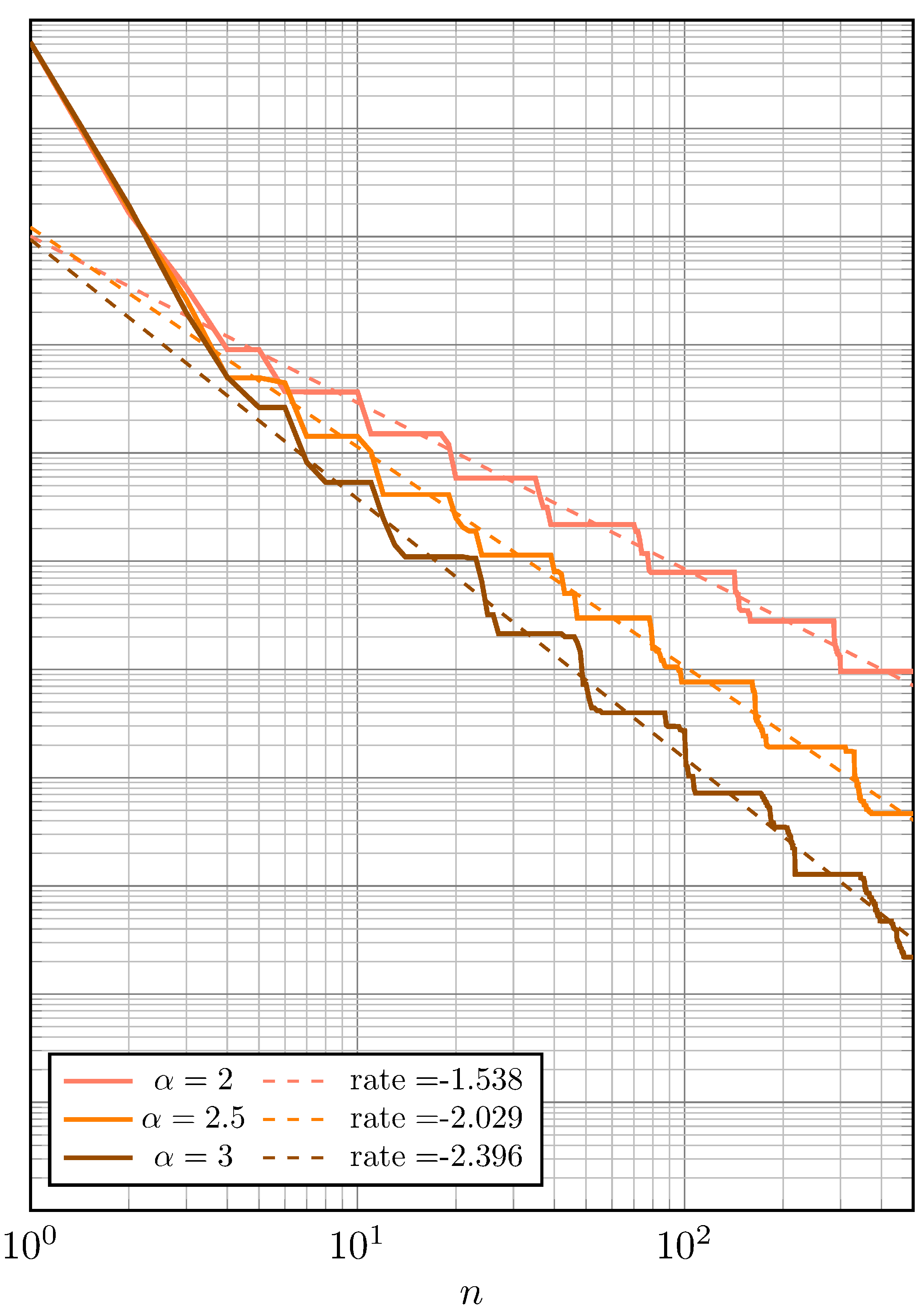}
		\caption{Mollifier mapping with wavelet expansion}
	\end{subfigure}\\
	\begin{subfigure}[t]{0.49\textwidth}
		\includegraphics[scale=0.074]{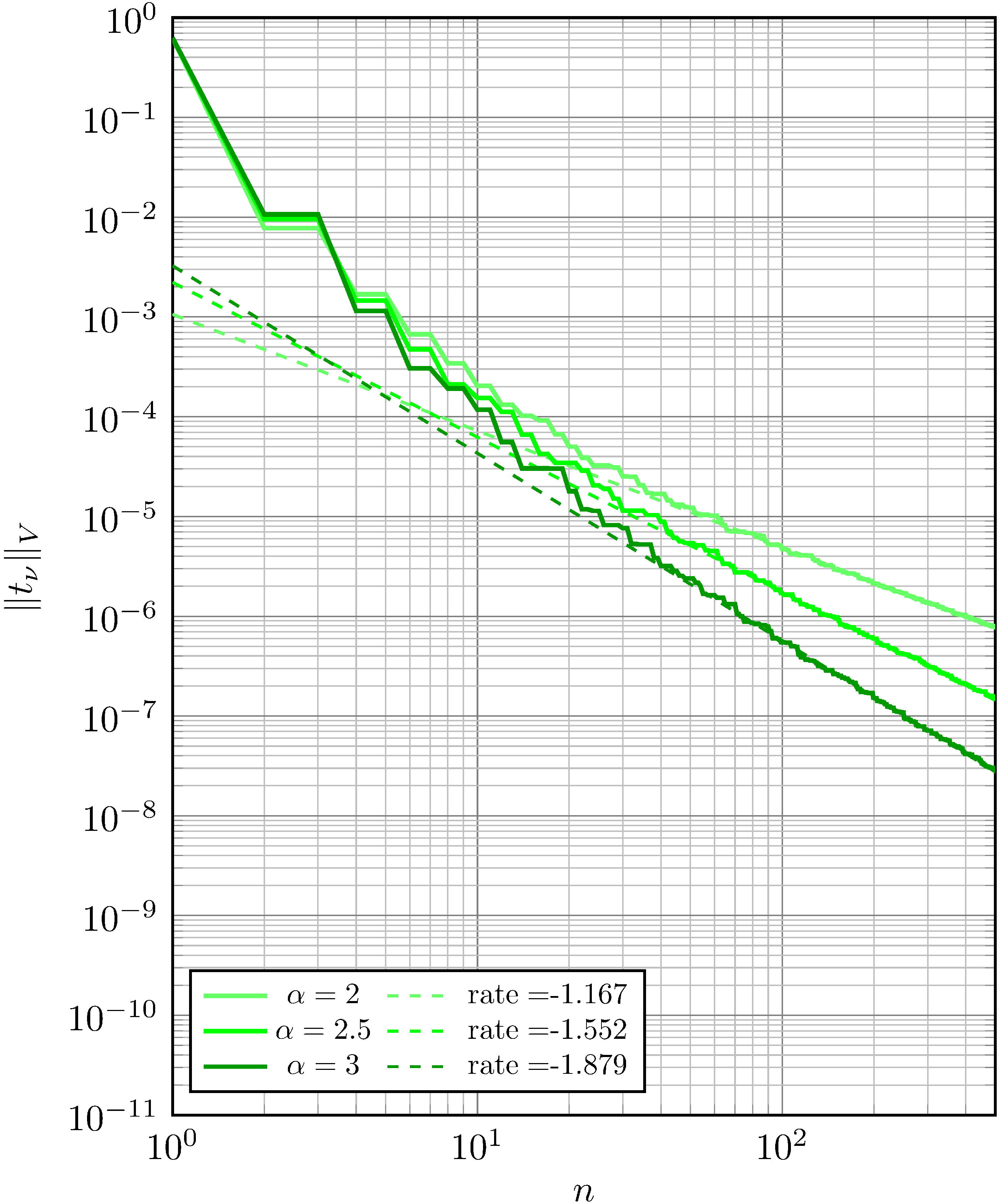}
		\caption{Harmonic mapping with Fourier expansion}
	\end{subfigure}
	~
	\begin{subfigure}[t]{0.45\textwidth}
		\includegraphics[scale=0.074]{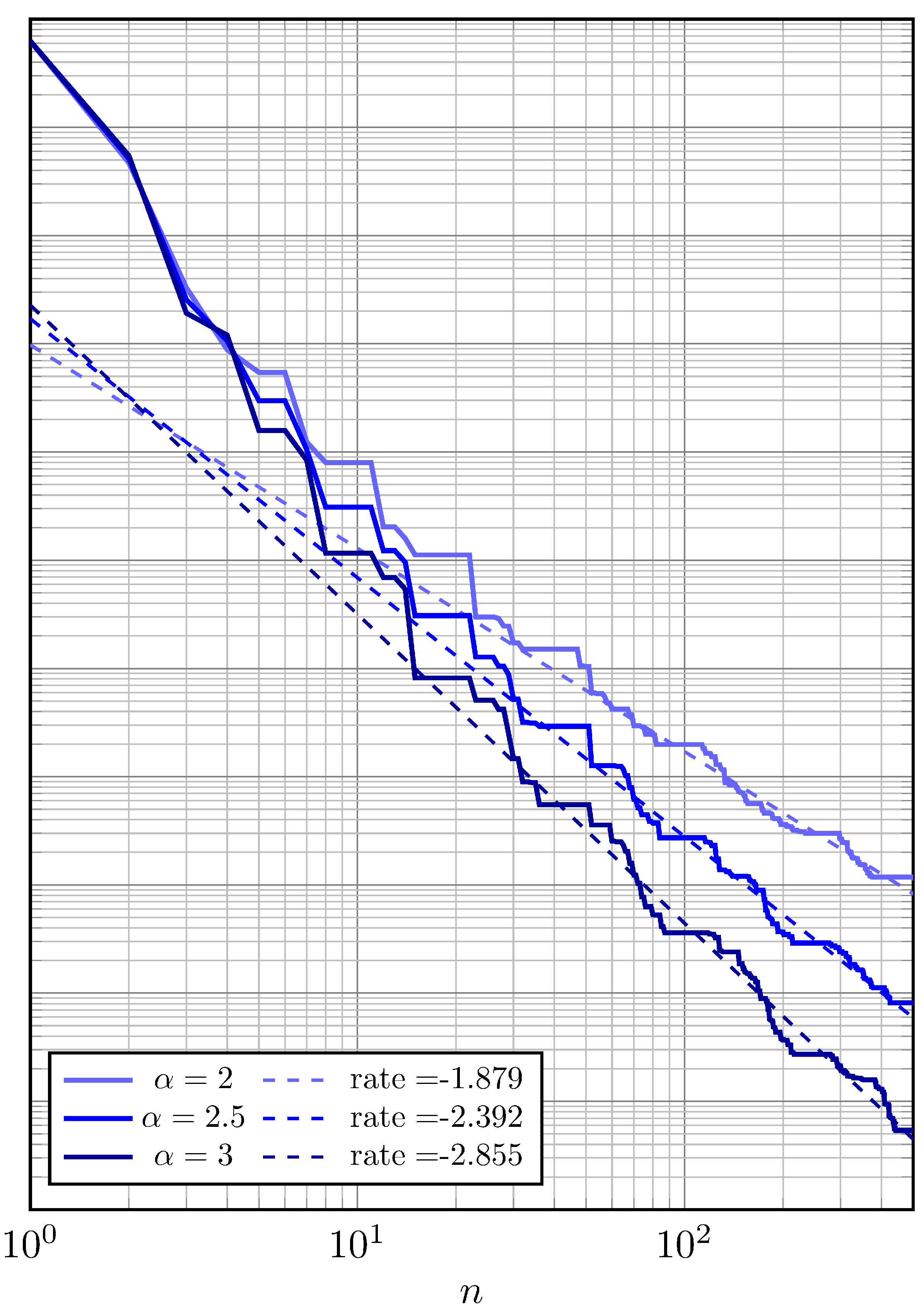}
		\caption{Harmonic mapping with wavelet expansion}
	\end{subfigure}
	\fi
	\caption{The effect of both the mapping and the expansion on the observed decay rate for $\vartheta=5\%$. We denote by $n$ the index in the decreasing rearrangement $(t^*_n)_{n \geq 1}$.}
	\label{fig:squareresults}
\end{figure}

\begin{table}
	\centering
	\begin{tabular}{l|cc|cc|cc|c}
			         & \multicolumn{2}{c|}{$h=4$}   	 		& \multicolumn{2}{c|}{$h=2$} 	& \multicolumn{2}{c|}{$h=1$} 			 			& $h=1/2$ 		 	\\ \hline
		$\vartheta$  & \multicolumn{1}{c|}{M} 	   & H   	  & \multicolumn{1}{c|}{M} 			& H   	   & \multicolumn{1}{c|}{M} 	 &H 		& H     		    	\\ \hline
		$4\%$        & \multicolumn{1}{l|}{2.512}  & 2.376    & \multicolumn{1}{l|}{2.457}  	& 2.629    & \multicolumn{1}{l|}{2.409}  & 2.867  	& 2.944     	      	\\
		$8\%$        & \multicolumn{1}{l|}{2.341}  & 2.496    & \multicolumn{1}{l|}{2.320}  	& 2.729    & \multicolumn{1}{l|}{2.408}  & 2.858  	& 2.877      	       	\\
		$16\%$       & \multicolumn{1}{l|}{2.431}  & 2.491    & \multicolumn{1}{l|}{2.421}  	& 2.629    & \multicolumn{1}{l|}{2.376}  & 2.700  	& 2.720      	      	\\
		$32\%$       & \multicolumn{1}{l|}{2.200}  & 2.369    & \multicolumn{1}{l|}{2.194}  	& 2.489    & \multicolumn{1}{l|}{2.193}  & 2.507  	& 2.507      	     
	\end{tabular}
	
	\caption{Decay of the $V$-norm of the Taylor coefficients of the PDE solution for the wavelet expansion with the mollifier (M) and harmonic (H) mapping, for different values of the maximal shape variation $\vartheta$ and fixed $\alpha=3$. Mesh sizes $h$ are relative to the mesh size required to adequately resolve the wavelets in Figure \ref{fig:threeresults} and \ref{fig:squareresults}.}
	\label{tab:shapevar}
\end{table}

\subsection{Illustrations of the theoretical rates for fixed maximal shape variations}
\label{sec:vartheta5percent} 

We verify the convergence rates obtained in Section \ref{sec:modelproblem} by computing the Taylor coefficients when using maximal shape variations $\vartheta=5\%$. We do this for both the mollifier-based and the harmonic mapping, where, in the first case, we use the mollifier introduced in Example \ref{def:linmollifier}. 

In Figure \ref{fig:threeresults}, we show the convergence graphs separated by value of the decay rate $\alpha$ for $\alpha\in \{2, 2.5, 3\}$. Previous results for Fourier and our results for wavelets predict convergence rates of 1, 1.5, and 2 for the Fourier expansion and 1.5, 2, 2.5 for the wavelet expansion, respectively. For the mollifier mapping, we recover the expected and predicted convergence rates for $\alpha\in\{2, 2.5\}$. For $\alpha=3$, we obtain convergence rates just short of the expected ones, analogously to \cite{Bachmayr2017a} . Although we do not obtain the full rates for all values of $\alpha$, we observe an improvement in the convergence rate of $\frac{1}{2}$ between the Fourier and wavelet expansions under the mollifier mapping, as expected by the theory. 

Next to the results for the mollifier mapping, the decay in the norm of the Taylor coefficients obtained with the harmonic mapping are also shown in Figure \ref{fig:threeresults}. We observe that, in the case of the Fourier expansion, the decay rates are very similar for both mappings, where the harmonic mapping introduces a significantly smaller constant when compared to the mollifier mapping. This difference can be attributed to the difference in $W^{1,\infty}$-norm between the mollifier and harmonic mappings.  

Contrary to the Fourier expansion, when considering the wavelet expansion for the radius, we observe different decay rates for the two domain mappings. While the wavelet expansion together with the mollifier mapping approaches the expected convergence rate, the harmonic mapping results in convergence rates exceeding the expected ones consistently by approximately $0.3$. This difference could be explained by the different locality behavior of the two mappings. Since the mollifier mappings are axially supported corresponding to the support of the wavelets, the support in the radial direction of the mollifier mapping is equal to the support of the mollifier $\chi$. This is different for the harmonic mapping, where not only the support in the angular direction shrinks as the support of the wavelet shrinks, but the support along the radial direction shrinks as well. Therefore, the harmonic mappings are more local in area when compared to the mollifier mapping. 

In order to investigate the effect of the decay parameter $\alpha$, we show the same results as in Figure \ref{fig:threeresults}, separated by mapping and expansion, in Figure \ref{fig:squareresults}. From this, we can see a clear difference between the mollifier and harmonic mappings. While the first two graphs corresponding to the mollifier mapping show almost no pre-asymptotic behavior, the equivalent graphs for the harmonic mapping show a pre-asymptotic behavior in the 50 largest coefficients. 

\subsection{Effect of maximal shape variations on the observed rates}
\label{sec:shapevar} 
Furthermore, we investigate the effect of the maximal shape variations on the obtained convergence rates. To this extent, we let the maximal shape variations be $\vartheta \in \{4\%, 8\%, 16\%, 32\% \}$, and we calculate the decay rate of the Taylor coefficients for each value of $\vartheta$, for both the linear mollifier and the harmonic mapping with the wavelet basis expansion. Table \ref{tab:shapevar} shows the results for different spatial resolutions, and in this subsection we focus on the column for $h=1$. 

Here, we observe that the difference between the decay rates for the mollifier and harmonic mapping persists for different values of $\vartheta$ and, for both of them, the decay rates degrade for larger values of $\vartheta$.  These observations are consistent with those in \cite{Bachmayr2017a} for the affine case. Moreover, the column for $h=1$ in Table \ref{tab:shapevar} shows us that, for small values of $\vartheta$, the decay rate when using the harmonic mapping approaches $3$, which is a significant improvement over the theoretical convergence rate of $2.5$. 

Additionally, we further investigate the structure of the active index set $\Lambda$ for different values of $\vartheta$. For small values of $\vartheta$, the active index set consists mainly of Taylor coefficients of small order ($|\mu|\leq 2$). On the other hand, for large values of $\vartheta$, the alternating greedy Taylor algorithm explores higher order coefficients earlier. 

\subsection{Effect of the space discretization on the observed rates}
\label{sec:convstudy} 
Next to the effect of the maximal shape variations $\vartheta$, the different columns in Table \ref{tab:shapevar} show the dependence of the rates on the spatial discretizations as well. In our experiments, $h=1$ is equivalent to the mesh refinement needed to resolve the wavelets in the $10^{\text{th}}$ layer on the boundary, as these are the smallest wavelets that enter the active set $\Lambda$ for the first 1000 iterations of the alternating greedy Taylor algorithm. Due to this choice, we do not expect any significant improvement for $h < 1$ and this is confirmed by the numerical experiments in Table \ref{tab:shapevar}. We have not performed the corresponding computations for the mollifier mapping, as the mesh would be prohibitively large due to the mesh size scaling decreasing linearly with the distance to the origin on the support of the mollifier, resulting in a factor 16 difference between the number of vertices in the mollifier and harmonic mesh.

Finally, we observe that the decay rates of the harmonic mapping decrease for larger values of $h$. This degradation is clearly noticeable for both small and large values of $\vartheta$, but the effect is of different magnitude: smaller values of $\vartheta$ imply larger degradation, and vice-versa. Contrary to the harmonic mapping, the decay rates for the mollifier mapping degrade, but by a smaller amount. This could also be due to the strong damping of the harmonic mapping inside the domain. Moreover, the degradation is not monotonic; if there is degradation, the decay rates become larger for some, and smaller for other values of $\vartheta$.

%% file: 6conclusion/6conclusion.tex
\section{Discussion and extensions}
\label{sec:conclusion}
In this work, we have shown that the convergence rates of Taylor coefficients of solutions to a certain class of elliptic PDEs can benefit greatly from using functions with localized support to represent function-valued parametric inputs. 
Namely, this can lead to a theoretical improvement of $\frac{1}{2}$ in the decay of the Taylor coefficients. 
Due to Stechkin's lemma, this  translates to an improvement of $\frac{1}{2}$ in the asymptotic convergence of the best-$N$-term truncation of the Taylor surrogate model. 

In application to parametric domains, addressed in Section \ref{sec:modelproblem}, we have discussed two different mapping approaches, an affine mollifier approach and a mapping defined through a harmonic extension problem. Analytically, the mollifier approach has an advantage over the harmonic mapping because we can establish pointwise bounds on the Jacobian matrix $\D\Phi$ by explicit calculations for star-shaped domains. This is in contrast to the harmonic extension, where we have shown similar bounds only for circular domains. However, intuitively, similar results for other classes of star-shaped domains are to be expected.

In contrast to the analytical preference for the mollifier mapping, the harmonic mapping has major numerical advantages over the first one. First of all, the harmonic mapping is such that the spectral norm of the Jacobian matrix $\D\Phi$  is relatively small compared to one for the mollifier mapping. Moreover, wavelets deep into the wavelet structure are such that the support of $\D\Phi_j$ is strongly localized near the domain boundary. Therefore, calculating the Taylor coefficients requires a fine mesh near the boundary of the computational domain to resolve the fine details in the small wavelets and can be kept rough on its interior to reduce the computational burden.

Instead of calculating the transformation by solving a partial differential equation for each dimension separately, other PDE-based approaches to build the mapping act on all coordinates simultaneously. One example of such a mapping is defined through the linear elasticity equations \cite{Cizmas2008,Dwight2009}, which one can interpret as elastically stretching the nominal onto the parameterized domain. Elasticity equations offer even more control over the smoothing properties of the transformation by using piecewise constant material properties. Although we have focused on star-shaped domains in this work, a significant advantage of PDE-based mapping constructions is their geometrical flexibility, including possible extensions to not star-shaped domains. 

We have not investigated it here, but, because of similarities in the PDEs, we expect our approach to work for the fourth-order equation $\Delta(a\Delta u)=f$ on a bounded Lipschitz domain, or for the parabolic equation $\partial_tu-\nabla \cdot(A\nabla u) = f$ set on $(0,T)\times D$. For these and similar equations, we refer to the discussion in \cite{Bachmayr2017a}. 

The theoretical convergence rate obtained in Theorem \ref{thm:lpsummability} is supported by numerical results for the parametric domain problem \eqref{eq:transformedvarform}. The numerical experiments show that both the choice for the basis expansion and the mapping approach majorly impact the observed convergence rate. Moreover, the crucial role of the mesh in the process showcases the possible benefits of parameter adapted meshes, where the mesh is refined adaptively, depending on the parameter. 

As mentioned in the introduction, our results can be valuable in uncertainty quantification, when, adopting a parametric approach, the parameter $\y$ is the image of independent, identically distributed random variables. However, we stress that, in that case, the choice of the functions in the expansion of the spatially-dependent input is a modeling choice: different expansions give rise to random fields with diverse statistical properties. 

A possible extension of our work is to consider the summability of Legendre coefficients, for which we expect the proof to be doable by building on the estimates in \cite{Bachmayr2017a,Dung2022}.

Finally, we want to mention the potential implications of our results. Since the summability properties of Taylor (and Legendre) expansions have often been used to show convergence rates of other polynomial surrogates, such as polynomial chaos- and collocation-based ones \cite{Beck2014,Beck2012,Chkifa2014}, their convergence rates may also benefit from expansions with localized support. For heuristic algorithms, this is not necessarily the case \cite{Ernst2021}. Next to this, they are also used to prove convergence of quadrature \cite{Dick2016,Schillings2013,Zech2020} for the computation of moments, approximation properties of neural networks \cite{Opschoor2022,Schwab2019}, and reduced basis and proper orthogonal decomposition \cite{Bachmayr2017a,Cohen2015}. The performance of these methods may also benefit from locality, but sound theoretical and numerical investigations for these frameworks are needed.

%% file: 7appendices/7appendices.tex
\begin{appendices}
\section{Combinatorial Lemma}
\label{ap:combilemma}

Whilst proving Lemma \ref{thm:Afbound}, we expand the derivatives in $|\A |_{2,2}$. When doing so, we need to apply the general Leibnitz rule several times, and this will results in sums of the type 
\begin{align*}
	\sum_{\nu^{(1)}+\cdots+\nu^{(N)}=\mu} \binom{\mu }{ \nu^{(1)},\cdots,\nu^{(N)}}\prod_{i=1}^N |\nu^{(i)}|!.
\end{align*}
To handle these sums, we introduce a combinatorial lemma.

\begin{lemma} \label{lem:multinomialsum}
	Let $\mu$ be a multi index with finite support and $N\in\mathbb{N}$. Now, we have:
	\begin{equation}
		\sum_{\nu^{(1)}+\cdots+\nu^{(N)}=\mu} \binom{\mu }{ \nu^{(1)},\cdots,\nu^{(N)}}\prod_{i=1}^N |\nu^{(i)}|! = \frac{(\abs{\mu}+(N-1))!}{(N-1)!}
	\end{equation}
	where $\binom{\mu }{ \nu^{(1)},\cdots,\nu^{(N)}}$ is the multinomial coefficient.
\end{lemma}

\begin{proof}
	Let $n=\argmax_{j \geq 1} \mu_j$. Now, let $S$ be a stack of cards in $n$ suits together with $N-1$ indistinghuisable jokers. For each suit $i$, we have $\mu_i$ cards per suit. A natural question to ask is: `In how many ways can deck $S$ be ordered?'
	
	First, we will calculate this directly. We can order $\abs{\mu}+(N-1)$ cards in $(\abs{\mu}+(N-1))!$ ways. To account for the $N-1$ identical jokers, we devide by $(N-1)!$ and we end up with 
	\begin{equation*}
		\# \text{orderings of }S = \frac{(\abs{\mu}+(N-1))!}{(N-1)!}.
	\end{equation*}
	
	Next, we calculate the number of orderings by summing over the number of cards per suit between the jokers. To this extent, we denote by $\nu_i^{(j)}$ the number of cards between jokers $j-1$ and $j$ of suit $i$. To account for all cards in the deck, we require $\sum_{j=1}^N \nu^{(j)}=\mu$. 
	 
	For each suit $i$, we can divide the $\mu_i$ cards into the $(\nu_i^{(j)})_{j=1}^N$ sections in $\binom{\mu_i }{ \nu_i^{(1)},\cdots,\nu_i^{(N)}}$ ways. Then we linearly order all cards before the first joker in $(\nu_1^{(1)}+\cdots + \nu_n^{(1)})!=|\nu^{(1)}|!$ ways, the cards between the first and the second joker in $|\nu^{(2)}|!$ ways, etcetera. Combining this, we can order all these card selections in $\prod_{j=1}^N |\nu^{(j)}|!$ ways. Combing this, we get 
	\begin{align*}
		&\binom{\mu_1 }{ \nu_1^{(1)},\cdots,\nu_1^{(N)}}\cdots \binom{\mu_n }{ \nu_n^{(1)},\cdots,\nu_n^{(N)}} \cdot |\nu^{(1)}|! \cdots |\nu^{(n)}|!\\
		&=\binom{\mu }{ \nu^{(1)},\cdots,\nu^{(N)}}\prod_{i=1}^N |\nu^{(i)}|!
	\end{align*}
	Finally, by summing over all possible $\nu^{(1)}+\cdots+\nu^{(N)} = \mu$ we get
	\begin{equation*}
		\# \text{orderings of }S = \sum_{\nu^{(1)}+\cdots+\nu^{(N)}=\mu} \binom{\mu }{ \nu^{(1)},\cdots,\nu^{(N)}}\prod_{i=1}^N |\nu^{(i)}|!
	\end{equation*}
	which shows the desired result.
\end{proof}

\section{Proof of Lemma \ref{lem:Dphibound}}
\label{pf:dphibound}
Here, we present the proof of Lemma \ref{lem:Dphibound}.
\begin{proof}
We rewrite equation \eqref{eq:mollifiertransf} to obtain
\begin{equation*}
	\Phi(\x;\bm{y})=\hat{\bm{x}} + \sum_{j \geq 1} \Phi_j(\x)y_j
\end{equation*}
with
\begin{equation*}
	\Phi_j(\x) = \chi\left(\x\right)  \psi_j(\theta(\x)) \frac{\x}{|\x|}.
\end{equation*}
Next, we take the derivative with respect to $\x$ and take the $|\cdot|_{2,2}$ norm:
\begin{align}
	 | \D \Phi_j(\x) |_{2,2} &= | \D \left[\chi\left(\x\right)  \psi_j(\theta(\x)) \frac{\x}{|\x|} \right] |_{2,2}\nonumber\\
	 &=|  \D\left[ \chi\left(\x\right) \frac{\x}{|\x|} \right] \psi_j(\theta(\x))  + 
	 \chi\left(\x\right) \frac{\x}{|\x|} \D\left[\theta(\x) \right]  \psi_j'(\theta(\x)) |_{2,2}\nonumber \\
	 &\leq| \psi_j(\theta(\x)) ||  A_1(\x) |_{2,2}  + 
	 |  \psi_j'(\theta(\x))  || A_2(\x)|_{2,2}, \label{eq:derbound}
\end{align}
with matrix-valued functions $A_1(\x)$ and $A_2(\x)$. First, we expand $A_1(\x)$:
\begin{align*}
	A_1(\x) &= \D\left[ \chi\left(\x\right) \frac{\x}{|\x|} \right] \\
	&=  \frac{1}{|\x|^3}\begin{bmatrix} |\x|^2\xx\chi_{x_1}(\x)+\chi\left(\x\right)  \xy^2  & |\x|^2\xx\chi_{x_2}(\x)-\chi\left(\x\right) \xx\xy \\ |\x|^2\xy\chi_{x_1}(\x) -\chi\left(\x\right) \xx\xy& |\x|^2\xy\chi_{x_2}(\x)+\chi\left(\x\right) \xx^2 \end{bmatrix}	,
\end{align*}
where the partial derivatives $\frac{\partial \chi}{\partial \xx}$ and $\frac{\partial \chi}{\partial \xy}$ are denoted by $\chi_{\xx}$ and  $\chi_{\xy}$ respectively. To obtain the $|\cdot|_{2,2}$-norm of $A_1(\x)$, we calculate the largest singular value $\sigma^+$ explicitly. First, we calculate the $A_1 A_1^\top$:
\begin{align*}
	A_1 A_1^\top&= \frac{1}{|\x|^6}\begin{bmatrix} |\x|^2\xx\chi_{x_1}(\x)+\chi\left(\x\right)  \xy^2  & |\x|^2\xx\chi_{x_2}(\x)-\chi\left(\x\right) \xx\xy \\ |\x|^2\xy\chi_{x_1}(\x) -\chi\left(\x\right) \xx\xy& |\x|^2\xy\chi_{\xy}(\x)+\chi\left(\x\right) \xx^2 \end{bmatrix} \cdot \\ &\qquad \qquad \begin{bmatrix} |\x|^2\xx\chi_{\xx}(\x)+\chi\left(\x\right)  \xy^2  &  |\x|^2\xy\chi_{\xx}(\x) -\chi\left(\x\right) \xx\xy\\|\x|^2\xx\chi_{\xy}(\x)-\chi\left(\x\right) \xx\xy  & |\x|^2\xy\chi_{\xy}(\x)+\chi\left(\x\right) \xx^2 \end{bmatrix}\\
	&=\frac{1}{|\x|^4}\left[  
\begin{array}{cc}
 \r^4 \chi_{\xx}(\x)^2+\xy^2 \chi(\x)^2 & \r^4 \chi_{\xy}(\x) \chi_{\xx}(\x)-\xx \xy \chi(\x)^2 \\
 \r^4 \chi_{\xy }(\x) \chi_{\xx}(\x)-\xx \xy \chi(\x)^2 & \r^4 \chi_{\xy}(\x)^2+\xx^2 \chi(\x)^2\\
\end{array}
\right] .
\end{align*}
To calculate the largest eigenvalue, we calcualte the trace:
\begin{align*}
	T&:= \tr(A_1 A_1^\top)
	= \frac{1}{|\x|^2}\left( \r^2  |\nabla_{\x} \chi|^2  +\chi(\x)^2   \right)
	\geq0
\end{align*}
and the determinant:
\begin{align*}
	D&:=\det(A_1 A_1^\top)  
	 =\frac{1}{\r^4} \chi(\x)^2 \left( \nabla_{\x}\chi(\x)\cdot \x \right) ^2 
	 \geq 0.
\end{align*}
Combining this, we obtain, for the largest singular value $\sigma_+$,
\begin{align*}
	\sigma_+^2 &= \frac{1}{2}\left(  T + \sqrt{T^2-4D} \right)\\
	&=\frac{1}{2}\cdot \\&\left(  \frac{1}{|\x|^2}\left( \r^2  |\nabla_{\x} \chi|^2  +\chi(\x)^2   \right) + \sqrt{\frac{1}{|\x|^4}\left( \r^2  |\nabla_{\x} \chi|^2  +\chi(\x)^2   \right)^2-\frac{4}{\r^4} \chi(\x)^2 \left( \nabla_{\x}\chi(\x)\cdot \x \right) ^2} \right),
\end{align*}
from which we conclude 
\begin{align}
	&\sigma_+ = \frac{1}{\sqrt{2}}\cdot\\&\sqrt{   \frac{1}{|\x|^2}\left( \r^2  |\nabla_{\x} \chi|^2  +\chi(\x)^2   \right) + \sqrt{\frac{1}{|\x|^4}\left( \r^2  |\nabla_{\x} \chi|^2  +\chi(\x)^2   \right)^2-\frac{4}{\r^4} \chi(\x)^2 \left( \nabla_{\x}\chi(\x)\cdot \x \right) ^2}  }.\label{eq:derbound1}
\end{align}

In a similar manner to our treatment of $A_1(\x)$, we expand $ A_2(\x)$:
\begin{align*}
	A_2(\x) &= \chi\left(\x\right)\frac{\x}{|\x|}\D\left[\theta(\x) \right] 
	=  \frac{\chi\left(\x\right)}{|\x|^3} \begin{bmatrix} -\xx\xy & \xx^2 \\ -\xy^2 &\xx\xy  \end{bmatrix} .
\end{align*}
For the determinant, we have 
\begin{align*}
	\det(A_2) &=   \frac{\chi\left(\x\right)}{|\x|^3}  \det \begin{bmatrix} -\xx\xy & \xx^2 \\ -\xy^2 &\xx\xy  \end{bmatrix} 
	  =0.
\end{align*}
Therefore, we calculate  the singular values of $B=\begin{bmatrix} -\xx\xy & \xx^2 \\ -\xy^2 &\xx\xy  \end{bmatrix}$ explicitly:
\begin{align*}
	0 &= \det\left( B^\top B -\lambda I_2\right)
	\det\left(\begin{bmatrix} -\xx\xy & -\xy^2 \\ \xx^2 &\xx\xy  \end{bmatrix}\begin{bmatrix} -\xx\xy & \xx^2 \\ -\xy^2 &\xx\xy  \end{bmatrix} -\lambda I_2 \right)
\end{align*}
which we expand to obtain
\begin{align*}	
	&\det \begin{bmatrix} \xx^2\xy^2 + \xy^4 - \lambda& -\xx^3\xy - \xx\xy^3 \\ -\xx^3\xy - \xx\xy^3 & \xx^4 + \xx^2\xy^2 -\lambda \end{bmatrix}  
	 = \lambda^2  - \lambda\left(\xx^2+\xy^2 \right)^2.
\end{align*}
From this, we arrive at two singular values, $\sigma_1 = 0$ and $\sigma_2=|\x|^2$ to get 
\begin{align}
|A_2(\x)|_{2,2} 
 &= \frac{\chi\left(\x\right)}{|\x|}.  \label{eq:derbound2}
\end{align}
Now, the result follows by combining equation \eqref{eq:derbound}, \eqref{eq:derbound1}, and  \eqref{eq:derbound2}.
\end{proof}

\section{Proof for bounds in Example \ref{col:Dphibound_bound}}
\label{ap:Dphibound_boundpf} 
\begin{proof}
	First, we expand and bound $\chi_1(\x)$:
	\begin{align}
		\chi_1(\x) &= \frac{1}{\sqrt{2}}\cdot\nonumber\\&\sqrt{   \frac{1}{|\x|^2}\left( \r^2  |\nabla_{\x} \chi|^2  +\chi(\x)^2   \right) + \sqrt{\frac{1}{|\x|^4}\left( \r^2  |\nabla_{\x} \chi|^2  +\chi(\x)^2   \right)^2-\frac{4}{\r^4} \chi(\x)^2 \left( \nabla_{\x}\chi(\x)\cdot \x \right) ^2}  }\nonumber\\
		&\leq \frac{1}{\sqrt{2}}\sqrt{   \frac{1}{|\x|^2}\left( \r^2  |\nabla_{\x} \chi|^2  +\chi(\x)^2   \right) + \sqrt{\frac{1}{|\x|^4}\left( \r^2  |\nabla_{\x} \chi|^2  +\chi(\x)^2   \right)^2}  }\nonumber\\
		&=  \sqrt{|\nabla_{\x} \chi|^2  +\frac{\chi(\x)^2}{\r^2}     }.\label{eq:singularvaluebound}
	\end{align}
	Next, we bound $|\nabla_{\x} \chi|^2$ with mollifier \eqref{eq:mollifier} for $\a \leq \r \leq r_0(\theta(\x))$:
	\begin{align}
		|\nabla_{\x} \chi|^2 
		&=\frac{\r^2\left(r_0(\theta(\x))-\a\right)^2 +\left( \r-\a    \right)^2  r_0'(\theta(\x)) ^2    }{\r^2\left(r_0(\theta(\x))-\a\right)^4}\nonumber\\
		&\leq \frac{1}{\left(r_0^- -\a\right)^2} + \frac{  r_0'(\theta(\x)) ^2    }{{r_0^-}^2\left(r_0^--\a\right)^2}\nonumber\\
		&= \frac{16}{9{r_0^-} ^2} \left( 1+ \frac{r_0'(\theta(\x)) ^2}{{r_0^-}^2}    \right)\label{eq:normgradsquaredbound}.
	\end{align}
	In a similar manner, we bound $\frac{\chi(\x)^2}{\r^2}$ in the $\a \leq \r \leq r_0(\theta(\x))$ case:
	\begin{align}
		\frac{\chi(\x)^2}{\r^2}&=\frac{\left(\r - \a\right)^2 }{\r^2\left( r_0(\theta(\x))-\a\right)^2}\nonumber\\
		&\leq\frac{16}{{9r_0^-}^2} \left(1 - \frac{r_0^-}{4r_0^+}\right)^2  \label{eq:chioverrbound}
	\end{align}
	Combinging equations \eqref{eq:singularvaluebound}, \eqref{eq:normgradsquaredbound}, and \eqref{eq:chioverrbound}, we obtain 
	\begin{align*}
		\chi_1(\x) &\leq \sqrt{\frac{16}{9{r_0^-} ^2} \left( 1+ \frac{r_0'(\theta(\x)) ^2}{{r_0^-}^2}    \right) + \frac{16}{{9r_0^-}^2} \left(1 - \frac{r_0^-}{4r_0^+}\right)^2 }\\
		&\leq \frac{4}{3{r_0^-}} \sqrt{ \left( 1+ \frac{L_{r_0} ^2}{{r_0^-}^2}    \right) + \left(1 - \frac{r_0^-}{4r_0^+}\right)^2 }\\
		&=\bar{\chi}_1,
	\end{align*}
	for $L_{r_0}$, the Lipschitz constant of $r_0(\theta)$.
	Finally, we close the proof by employing equation \eqref{eq:chioverrbound}, that is
	\begin{align*}
		\chi_2(\x)
		 \leq  \frac{4}{{3r_0^-}} \left(1 - \frac{r_0^-}{4r_0^+}\right) 
		 =\bar{\chi}_2.
	\end{align*}
\end{proof}

%
%
%

\end{appendices}